\input amstex
%Sekretariats-Stildatei
%
%AMSPPT.STY modifiziert und erweitert von C. Krattenthaler

\def\next{AMS-SEKR}\ifx\styname\next \endinput\fi
\catcode`\@=11
\def\styname{AMS-SEKR}
\def\styversion{2.0}
{\W@{}\W@{\styname.STY - Version \styversion}\W@{}}
\hyphenation{acad-e-my acad-e-mies af-ter-thought anom-aly anom-alies
an-ti-deriv-a-tive an-tin-o-my an-tin-o-mies apoth-e-o-ses apoth-e-o-sis
ap-pen-dix ar-che-typ-al as-sign-a-ble as-sist-ant-ship as-ymp-tot-ic
asyn-chro-nous at-trib-uted at-trib-ut-able bank-rupt bank-rupt-cy
bi-dif-fer-en-tial blue-print busier busiest cat-a-stroph-ic
cat-a-stroph-i-cally con-gress cross-hatched data-base de-fin-i-tive
de-riv-a-tive dis-trib-ute dri-ver dri-vers eco-nom-ics econ-o-mist
elit-ist equi-vari-ant ex-quis-ite ex-tra-or-di-nary flow-chart
for-mi-da-ble forth-right friv-o-lous ge-o-des-ic ge-o-det-ic geo-met-ric
griev-ance griev-ous griev-ous-ly hexa-dec-i-mal ho-lo-no-my ho-mo-thetic
ideals idio-syn-crasy in-fin-ite-ly in-fin-i-tes-i-mal ir-rev-o-ca-ble
key-stroke lam-en-ta-ble light-weight mal-a-prop-ism man-u-script
mar-gin-al meta-bol-ic me-tab-o-lism meta-lan-guage me-trop-o-lis
met-ro-pol-i-tan mi-nut-est mol-e-cule mono-chrome mono-pole mo-nop-oly
mono-spline mo-not-o-nous mul-ti-fac-eted mul-ti-plic-able non-euclid-ean
non-iso-mor-phic non-smooth par-a-digm par-a-bol-ic pa-rab-o-loid
pa-ram-e-trize para-mount pen-ta-gon phe-nom-e-non post-script pre-am-ble
pro-ce-dur-al pro-hib-i-tive pro-hib-i-tive-ly pseu-do-dif-fer-en-tial
pseu-do-fi-nite pseu-do-nym qua-drat-ics quad-ra-ture qua-si-smooth
qua-si-sta-tion-ary qua-si-tri-an-gu-lar quin-tes-sence quin-tes-sen-tial
re-arrange-ment rec-tan-gle ret-ri-bu-tion retro-fit retro-fit-ted
right-eous right-eous-ness ro-bot ro-bot-ics sched-ul-ing se-mes-ter
semi-def-i-nite semi-ho-mo-thet-ic set-up se-vere-ly side-step sov-er-eign
spe-cious spher-oid spher-oid-al star-tling star-tling-ly
sta-tis-tics sto-chas-tic straight-est strange-ness strat-a-gem strong-hold
sum-ma-ble symp-to-matic syn-chro-nous topo-graph-i-cal tra-vers-a-ble
tra-ver-sal tra-ver-sals treach-ery turn-around un-at-tached un-err-ing-ly
white-space wide-spread wing-spread wretch-ed wretch-ed-ly Brown-ian
Eng-lish Euler-ian Feb-ru-ary Gauss-ian Grothen-dieck Hamil-ton-ian
Her-mit-ian Jan-u-ary Japan-ese Kor-te-weg Le-gendre Lip-schitz
Lip-schitz-ian Mar-kov-ian Noe-ther-ian No-vem-ber Rie-mann-ian
Schwarz-schild Sep-tem-ber
%Zustze (Sind auch in besp.exc einzutragen!):
form per-iods Uni-ver-si-ty cri-ti-sism for-ma-lism}
\Invalid@\nofrills
\Invalid@\usualspace
\newif\ifnofrills@
\def\nofrills@#1#2{\relaxnext@
  \DN@{\ifx\next\nofrills
    \nofrills@true\let#2\relax\DN@\nofrills{\nextii@}%
  \else
    \nofrills@false\def#2{#1}\let\next@\nextii@\fi
\next@}}
\def\usualspace@#1{\ifnofrills@\def\usualspace{#1}\fi}
\def\addto#1#2{\csname \expandafter\eat@\string#1@\endcsname
  \expandafter{\the\csname \expandafter\eat@\string#1@\endcsname#2}}
\newdimen\bigsize@
\def\big@#1#2{{\hbox{$\left#2\vcenter to#1\bigsize@{}%
  \right.\nulldelimiterspace\z@\m@th$}}}
\def\big{\big@\@ne}
\def\Big{\big@{1.5}}
\def\bigg{\big@\tw@}
\def\Bigg{\big@{2.5}}
\def\raggedcenter@{\leftskip\z@ plus.4\hsize \rightskip\leftskip
 \parfillskip\z@ \parindent\z@ \spaceskip.3333em \xspaceskip.5em
 \pretolerance9999\tolerance9999 \exhyphenpenalty\@M
 \hyphenpenalty\@M \let\\\linebreak}
\def\upperspecialchars{\def\ss{SS}\let\i=I\let\j=J\let\ae\AE\let\oe\OE
  \let\o\O\let\aa\AA\let\l\L}
\def\uppercasetext@#1{%
  {\spaceskip1.2\fontdimen2\the\font plus1.2\fontdimen3\the\font
   \upperspecialchars\uctext@#1$\m@th\aftergroup\eat@$}}
\def\uctext@#1$#2${\endash@#1-\endash@$#2$\uctext@}
\def\endash@#1-#2\endash@{\uppercase{#1}\if\notempty{#2}--\endash@#2\endash@\fi}
\def\runaway@#1{\DN@{#1}\ifx\envir@\next@
  \Err@{You seem to have a missing or misspelled \string\end#1 ...}%
  \let\envir@\empty\fi}
\newif\iftemp@
\def\notempty#1{TT\fi\def\test@{#1}\ifx\test@\empty\temp@false
  \else\temp@true\fi \iftemp@}
\font@\tensmc=cmcsc10
\font@\sevenex=cmex7
\font@\sevenit=cmti7
\font@\eightrm=cmr8 % preloaded in plain.tex
\font@\sixrm=cmr6 % preloaded in plain.tex
\font@\eighti=cmmi8     \skewchar\eighti='177 % preloaded
\font@\sixi=cmmi6       \skewchar\sixi='177   % preloaded
\font@\eightsy=cmsy8    \skewchar\eightsy='60 % preloaded
\font@\sixsy=cmsy6      \skewchar\sixsy='60   % preloaded
\font@\eightex=cmex8
\font@\eightbf=cmbx8 % preloaded in plain.tex
\font@\sixbf=cmbx6   % preloaded in plain.tex
\font@\eightit=cmti8 % preloaded in plain.tex
\font@\eightsl=cmsl8 % preloaded in plain.tex
\font@\eightsmc=cmcsc8
\font@\eighttt=cmtt8 % preloaded in plain.tex
%\font@\ninerm=cmr9
%\font@\ninei=cmmi9    \skewchar\ninei='177
%\font@\ninesy=cmsy9   \skewchar\ninesy='60
%\font@\nineex=cmex9
%\font@\ninebf=cmbx9
%\font@\nineit=cmti9
%\font@\ninesl=cmsl9
%\font@\ninesmc=cmcsc9
%\font@\ninemsa=msam9
%\font@\ninemsb=msbm9
%\font@\nineeufm=eufm9

%Ergnzung des fetten Small-Capitals-Fonts:
%\font@\eightbsmc=cmbcsc10 scaled 833
%\font@\tenbsmc=cmbcsc10
%\font@\twelvebsmc=cmbcsc10 scaled \magstep1
%\font@\fourteenbsmc=cmbcsc10 scaled \magstep2
%\font@\seventeenbsmc=cmbcsc10 scaled \magstep3
%\font@\twentybsmc=cmbcsc10 scaled \magstep4

\loadmsam
\loadmsbm
\loadeufm
\UseAMSsymbols
\newtoks\tenpoint@
\def\tenpoint{\normalbaselineskip12\p@
 \abovedisplayskip12\p@ plus3\p@ minus9\p@
 \belowdisplayskip\abovedisplayskip
 \abovedisplayshortskip\z@ plus3\p@
 \belowdisplayshortskip7\p@ plus3\p@ minus4\p@
 \textonlyfont@\rm\tenrm \textonlyfont@\it\tenit
 \textonlyfont@\sl\tensl \textonlyfont@\bf\tenbf
 \textonlyfont@\smc\tensmc \textonlyfont@\tt\tentt
%Ergnzung des fetten Small-Capitals-Fonts:
 \textonlyfont@\bsmc\tenbsmc
 \ifsyntax@ \def\big##1{{\hbox{$\left##1\right.$}}}%
  \let\Big\big \let\bigg\big \let\Bigg\big
 \else
  \textfont\z@=\tenrm  \scriptfont\z@=\sevenrm  \scriptscriptfont\z@=\fiverm
  \textfont\@ne=\teni  \scriptfont\@ne=\seveni  \scriptscriptfont\@ne=\fivei
  \textfont\tw@=\tensy \scriptfont\tw@=\sevensy \scriptscriptfont\tw@=\fivesy
  \textfont\thr@@=\tenex \scriptfont\thr@@=\sevenex
        \scriptscriptfont\thr@@=\sevenex
  \textfont\itfam=\tenit \scriptfont\itfam=\sevenit
        \scriptscriptfont\itfam=\sevenit
  \textfont\bffam=\tenbf \scriptfont\bffam=\sevenbf
        \scriptscriptfont\bffam=\fivebf
  \setbox\strutbox\hbox{\vrule height8.5\p@ depth3.5\p@ width\z@}%
  \setbox\strutbox@\hbox{\lower.5\normallineskiplimit\vbox{%
        \kern-\normallineskiplimit\copy\strutbox}}%
 \setbox\z@\vbox{\hbox{$($}\kern\z@}\bigsize@=1.2\ht\z@
 \fi
 \normalbaselines\rm\ex@.2326ex\jot3\ex@\the\tenpoint@}
\newtoks\eightpoint@
\def\eightpoint{\normalbaselineskip10\p@
 \abovedisplayskip10\p@ plus2.4\p@ minus7.2\p@
 \belowdisplayskip\abovedisplayskip
 \abovedisplayshortskip\z@ plus2.4\p@
 \belowdisplayshortskip5.6\p@ plus2.4\p@ minus3.2\p@
 \textonlyfont@\rm\eightrm \textonlyfont@\it\eightit
 \textonlyfont@\sl\eightsl \textonlyfont@\bf\eightbf
 \textonlyfont@\smc\eightsmc \textonlyfont@\tt\eighttt
%Ergnzung des fetten Small-Capitals-Fonts:
 \textonlyfont@\bsmc\eightbsmc
 \ifsyntax@\def\big##1{{\hbox{$\left##1\right.$}}}%
  \let\Big\big \let\bigg\big \let\Bigg\big
 \else
  \textfont\z@=\eightrm \scriptfont\z@=\sixrm \scriptscriptfont\z@=\fiverm
  \textfont\@ne=\eighti \scriptfont\@ne=\sixi \scriptscriptfont\@ne=\fivei
  \textfont\tw@=\eightsy \scriptfont\tw@=\sixsy \scriptscriptfont\tw@=\fivesy
  \textfont\thr@@=\eightex \scriptfont\thr@@=\sevenex
   \scriptscriptfont\thr@@=\sevenex
  \textfont\itfam=\eightit \scriptfont\itfam=\sevenit
   \scriptscriptfont\itfam=\sevenit
  \textfont\bffam=\eightbf \scriptfont\bffam=\sixbf
   \scriptscriptfont\bffam=\fivebf
 \setbox\strutbox\hbox{\vrule height7\p@ depth3\p@ width\z@}%
 \setbox\strutbox@\hbox{\raise.5\normallineskiplimit\vbox{%
   \kern-\normallineskiplimit\copy\strutbox}}%
 \setbox\z@\vbox{\hbox{$($}\kern\z@}\bigsize@=1.2\ht\z@
 \fi
 \normalbaselines\eightrm\ex@.2326ex\jot3\ex@\the\eightpoint@}

%Definition von 12-point, 14-point und 17-point Fonts
\font@\twelverm=cmr10 scaled\magstep1
\font@\twelveit=cmti10 scaled\magstep1
\font@\twelvesl=cmsl10 scaled\magstep1
\font@\twelvesmc=cmcsc10 scaled\magstep1
\font@\twelvett=cmtt10 scaled\magstep1
\font@\twelvebf=cmbx10 scaled\magstep1
\font@\twelvei=cmmi10 scaled\magstep1
\font@\twelvesy=cmsy10 scaled\magstep1
\font@\twelveex=cmex10 scaled\magstep1
\font@\twelvemsa=msam10 scaled\magstep1
\font@\twelveeufm=eufm10 scaled\magstep1
\font@\twelvemsb=msbm10 scaled\magstep1
\newtoks\twelvepoint@
\def\twelvepoint{\normalbaselineskip15\p@
 \abovedisplayskip15\p@ plus3.6\p@ minus10.8\p@
 \belowdisplayskip\abovedisplayskip
 \abovedisplayshortskip\z@ plus3.6\p@
 \belowdisplayshortskip8.4\p@ plus3.6\p@ minus4.8\p@
 \textonlyfont@\rm\twelverm \textonlyfont@\it\twelveit
 \textonlyfont@\sl\twelvesl \textonlyfont@\bf\twelvebf
 \textonlyfont@\smc\twelvesmc \textonlyfont@\tt\twelvett
%Ergnzung des fetten Small-Capitals-Fonts:
 \textonlyfont@\bsmc\twelvebsmc
 \ifsyntax@ \def\big##1{{\hbox{$\left##1\right.$}}}%
  \let\Big\big \let\bigg\big \let\Bigg\big
 \else
  \textfont\z@=\twelverm  \scriptfont\z@=\tenrm  \scriptscriptfont\z@=\sevenrm
  \textfont\@ne=\twelvei  \scriptfont\@ne=\teni  \scriptscriptfont\@ne=\seveni
  \textfont\tw@=\twelvesy \scriptfont\tw@=\tensy \scriptscriptfont\tw@=\sevensy
  \textfont\thr@@=\twelveex \scriptfont\thr@@=\tenex
        \scriptscriptfont\thr@@=\tenex
  \textfont\itfam=\twelveit \scriptfont\itfam=\tenit
        \scriptscriptfont\itfam=\tenit
  \textfont\bffam=\twelvebf \scriptfont\bffam=\tenbf
        \scriptscriptfont\bffam=\sevenbf
  \setbox\strutbox\hbox{\vrule height10.2\p@ depth4.2\p@ width\z@}%
  \setbox\strutbox@\hbox{\lower.6\normallineskiplimit\vbox{%
        \kern-\normallineskiplimit\copy\strutbox}}%
 \setbox\z@\vbox{\hbox{$($}\kern\z@}\bigsize@=1.4\ht\z@
 \fi
 \normalbaselines\rm\ex@.2326ex\jot3.6\ex@\the\twelvepoint@}

\def\headfonts{\twelvepoint\bf}

\font@\fourteenrm=cmr10 scaled\magstep2
\font@\fourteenit=cmti10 scaled\magstep2
\font@\fourteensl=cmsl10 scaled\magstep2
\font@\fourteensmc=cmcsc10 scaled\magstep2
\font@\fourteentt=cmtt10 scaled\magstep2
\font@\fourteenbf=cmbx10 scaled\magstep2
\font@\fourteeni=cmmi10 scaled\magstep2
\font@\fourteensy=cmsy10 scaled\magstep2
\font@\fourteenex=cmex10 scaled\magstep2
\font@\fourteenmsa=msam10 scaled\magstep2
\font@\fourteeneufm=eufm10 scaled\magstep2
\font@\fourteenmsb=msbm10 scaled\magstep2
\newtoks\fourteenpoint@
\def\fourteenpoint{\normalbaselineskip15\p@
 \abovedisplayskip18\p@ plus4.3\p@ minus12.9\p@
 \belowdisplayskip\abovedisplayskip
 \abovedisplayshortskip\z@ plus4.3\p@
 \belowdisplayshortskip10.1\p@ plus4.3\p@ minus5.8\p@
 \textonlyfont@\rm\fourteenrm \textonlyfont@\it\fourteenit
 \textonlyfont@\sl\fourteensl \textonlyfont@\bf\fourteenbf
 \textonlyfont@\smc\fourteensmc \textonlyfont@\tt\fourteentt
%Ergnzung des fetten Small-Capitals-Fonts:
 \textonlyfont@\bsmc\fourteenbsmc
 \ifsyntax@ \def\big##1{{\hbox{$\left##1\right.$}}}%
  \let\Big\big \let\bigg\big \let\Bigg\big
 \else
  \textfont\z@=\fourteenrm  \scriptfont\z@=\twelverm  \scriptscriptfont\z@=\tenrm
  \textfont\@ne=\fourteeni  \scriptfont\@ne=\twelvei  \scriptscriptfont\@ne=\teni
  \textfont\tw@=\fourteensy \scriptfont\tw@=\twelvesy \scriptscriptfont\tw@=\tensy
  \textfont\thr@@=\fourteenex \scriptfont\thr@@=\twelveex
        \scriptscriptfont\thr@@=\twelveex
  \textfont\itfam=\fourteenit \scriptfont\itfam=\twelveit
        \scriptscriptfont\itfam=\twelveit
  \textfont\bffam=\fourteenbf \scriptfont\bffam=\twelvebf
        \scriptscriptfont\bffam=\tenbf
  \setbox\strutbox\hbox{\vrule height12.2\p@ depth5\p@ width\z@}%
  \setbox\strutbox@\hbox{\lower.72\normallineskiplimit\vbox{%
        \kern-\normallineskiplimit\copy\strutbox}}%
 \setbox\z@\vbox{\hbox{$($}\kern\z@}\bigsize@=1.7\ht\z@
 \fi
 \normalbaselines\rm\ex@.2326ex\jot4.3\ex@\the\fourteenpoint@}

\def\chapheadfonts{\fourteenpoint\bf}

\font@\seventeenrm=cmr10 scaled\magstep3
\font@\seventeenit=cmti10 scaled\magstep3
\font@\seventeensl=cmsl10 scaled\magstep3
\font@\seventeensmc=cmcsc10 scaled\magstep3
\font@\seventeentt=cmtt10 scaled\magstep3
\font@\seventeenbf=cmbx10 scaled\magstep3
\font@\seventeeni=cmmi10 scaled\magstep3
\font@\seventeensy=cmsy10 scaled\magstep3
\font@\seventeenex=cmex10 scaled\magstep3
\font@\seventeenmsa=msam10 scaled\magstep3
\font@\seventeeneufm=eufm10 scaled\magstep3
\font@\seventeenmsb=msbm10 scaled\magstep3
\newtoks\seventeenpoint@
\def\seventeenpoint{\normalbaselineskip18\p@
 \abovedisplayskip21.6\p@ plus5.2\p@ minus15.4\p@
 \belowdisplayskip\abovedisplayskip
 \abovedisplayshortskip\z@ plus5.2\p@
 \belowdisplayshortskip12.1\p@ plus5.2\p@ minus7\p@
 \textonlyfont@\rm\seventeenrm \textonlyfont@\it\seventeenit
 \textonlyfont@\sl\seventeensl \textonlyfont@\bf\seventeenbf
 \textonlyfont@\smc\seventeensmc \textonlyfont@\tt\seventeentt
%Ergnzung des fetten Small-Capitals-Fonts:
 \textonlyfont@\bsmc\seventeenbsmc
 \ifsyntax@ \def\big##1{{\hbox{$\left##1\right.$}}}%
  \let\Big\big \let\bigg\big \let\Bigg\big
 \else
  \textfont\z@=\seventeenrm  \scriptfont\z@=\fourteenrm  \scriptscriptfont\z@=\twelverm
  \textfont\@ne=\seventeeni  \scriptfont\@ne=\fourteeni  \scriptscriptfont\@ne=\twelvei
  \textfont\tw@=\seventeensy \scriptfont\tw@=\fourteensy \scriptscriptfont\tw@=\twelvesy
  \textfont\thr@@=\seventeenex \scriptfont\thr@@=\fourteenex
        \scriptscriptfont\thr@@=\fourteenex
  \textfont\itfam=\seventeenit \scriptfont\itfam=\fourteenit
        \scriptscriptfont\itfam=\fourteenit
  \textfont\bffam=\seventeenbf \scriptfont\bffam=\fourteenbf
        \scriptscriptfont\bffam=\twelvebf
  \setbox\strutbox\hbox{\vrule height14.6\p@ depth6\p@ width\z@}%
  \setbox\strutbox@\hbox{\lower.86\normallineskiplimit\vbox{%
        \kern-\normallineskiplimit\copy\strutbox}}%
 \setbox\z@\vbox{\hbox{$($}\kern\z@}\bigsize@=2\ht\z@
 \fi
 \normalbaselines\rm\ex@.2326ex\jot5.2\ex@\the\seventeenpoint@}

\font@\rrrrrm=cmr10 scaled\magstep4
\font@\bigtitlefont=cmbx10 scaled\magstep4

\parindent1pc
\normallineskiplimit\p@
\newdimen\indenti \indenti=2pc
\def\pageheight#1{\vsize#1}
\def\pagewidth#1{\hsize#1%
   \captionwidth@\hsize \advance\captionwidth@-2\indenti}
\pagewidth{30pc} \pageheight{47pc}
\def\topmatter{%
 \ifx\undefined\msafam
 \else\font@\eightmsa=msam8 \font@\sixmsa=msam6
   \ifsyntax@\else \addto\tenpoint{\textfont\msafam=\tenmsa
              \scriptfont\msafam=\sevenmsa \scriptscriptfont\msafam=\fivemsa}%
     \addto\eightpoint{\textfont\msafam=\eightmsa \scriptfont\msafam=\sixmsa
              \scriptscriptfont\msafam=\fivemsa}%
   \fi
 \fi
 \ifx\undefined\msbfam
 \else\font@\eightmsb=msbm8 \font@\sixmsb=msbm6
   \ifsyntax@\else \addto\tenpoint{\textfont\msbfam=\tenmsb
         \scriptfont\msbfam=\sevenmsb \scriptscriptfont\msbfam=\fivemsb}%
     \addto\eightpoint{\textfont\msbfam=\eightmsb \scriptfont\msbfam=\sixmsb
         \scriptscriptfont\msbfam=\fivemsb}%
   \fi
 \fi
 \ifx\undefined\eufmfam
 \else \font@\eighteufm=eufm8 \font@\sixeufm=eufm6
   \ifsyntax@\else \addto\tenpoint{\textfont\eufmfam=\teneufm
       \scriptfont\eufmfam=\seveneufm \scriptscriptfont\eufmfam=\fiveeufm}%
     \addto\eightpoint{\textfont\eufmfam=\eighteufm
       \scriptfont\eufmfam=\sixeufm \scriptscriptfont\eufmfam=\fiveeufm}%
   \fi
 \fi
 \ifx\undefined\eufbfam
 \else \font@\eighteufb=eufb8 \font@\sixeufb=eufb6
   \ifsyntax@\else \addto\tenpoint{\textfont\eufbfam=\teneufb
      \scriptfont\eufbfam=\seveneufb \scriptscriptfont\eufbfam=\fiveeufb}%
    \addto\eightpoint{\textfont\eufbfam=\eighteufb
      \scriptfont\eufbfam=\sixeufb \scriptscriptfont\eufbfam=\fiveeufb}%
   \fi
 \fi
 \ifx\undefined\eusmfam
 \else \font@\eighteusm=eusm8 \font@\sixeusm=eusm6
   \ifsyntax@\else \addto\tenpoint{\textfont\eusmfam=\teneusm
       \scriptfont\eusmfam=\seveneusm \scriptscriptfont\eusmfam=\fiveeusm}%
     \addto\eightpoint{\textfont\eusmfam=\eighteusm
       \scriptfont\eusmfam=\sixeusm \scriptscriptfont\eusmfam=\fiveeusm}%
   \fi
 \fi
 \ifx\undefined\eusbfam
 \else \font@\eighteusb=eusb8 \font@\sixeusb=eusb6
   \ifsyntax@\else \addto\tenpoint{\textfont\eusbfam=\teneusb
       \scriptfont\eusbfam=\seveneusb \scriptscriptfont\eusbfam=\fiveeusb}%
     \addto\eightpoint{\textfont\eusbfam=\eighteusb
       \scriptfont\eusbfam=\sixeusb \scriptscriptfont\eusbfam=\fiveeusb}%
   \fi
 \fi
 \ifx\undefined\eurmfam
 \else \font@\eighteurm=eurm8 \font@\sixeurm=eurm6
   \ifsyntax@\else \addto\tenpoint{\textfont\eurmfam=\teneurm
       \scriptfont\eurmfam=\seveneurm \scriptscriptfont\eurmfam=\fiveeurm}%
     \addto\eightpoint{\textfont\eurmfam=\eighteurm
       \scriptfont\eurmfam=\sixeurm \scriptscriptfont\eurmfam=\fiveeurm}%
   \fi
 \fi
 \ifx\undefined\eurbfam
 \else \font@\eighteurb=eurb8 \font@\sixeurb=eurb6
   \ifsyntax@\else \addto\tenpoint{\textfont\eurbfam=\teneurb
       \scriptfont\eurbfam=\seveneurb \scriptscriptfont\eurbfam=\fiveeurb}%
    \addto\eightpoint{\textfont\eurbfam=\eighteurb
       \scriptfont\eurbfam=\sixeurb \scriptscriptfont\eurbfam=\fiveeurb}%
   \fi
 \fi
 \ifx\undefined\cmmibfam
 \else \font@\eightcmmib=cmmib8 \font@\sixcmmib=cmmib6
   \ifsyntax@\else \addto\tenpoint{\textfont\cmmibfam=\tencmmib
       \scriptfont\cmmibfam=\sevencmmib \scriptscriptfont\cmmibfam=\fivecmmib}%
    \addto\eightpoint{\textfont\cmmibfam=\eightcmmib
       \scriptfont\cmmibfam=\sixcmmib \scriptscriptfont\cmmibfam=\fivecmmib}%
   \fi
 \fi
 \ifx\undefined\cmbsyfam
 \else \font@\eightcmbsy=cmbsy8 \font@\sixcmbsy=cmbsy6
   \ifsyntax@\else \addto\tenpoint{\textfont\cmbsyfam=\tencmbsy
      \scriptfont\cmbsyfam=\sevencmbsy \scriptscriptfont\cmbsyfam=\fivecmbsy}%
    \addto\eightpoint{\textfont\cmbsyfam=\eightcmbsy
      \scriptfont\cmbsyfam=\sixcmbsy \scriptscriptfont\cmbsyfam=\fivecmbsy}%
   \fi
 \fi
 \let\topmatter\relax}
\def\chapterno@{\uppercase\expandafter{\romannumeral\chaptercount@}}
\newcount\chaptercount@
\def\chapter{\nofrills@{\afterassignment\chapterno@
                        CHAPTER \global\chaptercount@=}\chapter@
 \DNii@##1{\leavevmode\hskip-\leftskip
   \rlap{\vbox to\z@{\vss\centerline{\eightpoint
   \chapter@##1\unskip}\baselineskip2pc\null}}\hskip\leftskip
   \nofrills@false}%
 \FN@\next@}
\newbox\titlebox@

%Uppercase ist abgestellt.
\def\title{\nofrills@{\relax}\title@%
 \DNii@##1\endtitle{\global\setbox\titlebox@\vtop{\tenpoint\bf
 \raggedcenter@\ignorespaces
 \baselineskip1.3\baselineskip\title@{##1}\endgraf}%
 \ifmonograph@ \edef\next{\the\leftheadtoks}\ifx\next\empty
    \leftheadtext{##1}\fi
 \fi
 \edef\next{\the\rightheadtoks}\ifx\next\empty \rightheadtext{##1}\fi
 }\FN@\next@}
\newbox\authorbox@
\def\author#1\endauthor{\global\setbox\authorbox@
 \vbox{\tenpoint\smc\raggedcenter@\ignorespaces
 #1\endgraf}\relaxnext@ \edef\next{\the\leftheadtoks}%
 \ifx\next\empty\leftheadtext{#1}\fi}
\newbox\affilbox@
\def\affil#1\endaffil{\global\setbox\affilbox@
 \vbox{\tenpoint\raggedcenter@\ignorespaces#1\endgraf}}
\newcount\addresscount@
\addresscount@\z@
\def\address#1\endaddress{\global\advance\addresscount@\@ne
  \expandafter\gdef\csname address\number\addresscount@\endcsname
  {\vskip12\p@ minus6\p@\noindent\eightpoint\smc\ignorespaces#1\par}}
\def\email{\nofrills@{\eightpoint{\it E-mail\/}:\enspace}\email@
  \DNii@##1\endemail{%
  \expandafter\gdef\csname email\number\addresscount@\endcsname
  {\def\usualspace{{\it\enspace}}\smallskip\noindent\eightpoint\email@
  \ignorespaces##1\par}}%
 \FN@\next@}
\def\thedate@{}
\def\date#1\enddate{\gdef\thedate@{\tenpoint\ignorespaces#1\unskip}}
\def\thethanks@{}
\def\thanks#1\endthanks{\gdef\thethanks@{\eightpoint\ignorespaces#1.\unskip}}
\def\thekeywords@{}
\def\keywords{\nofrills@{{\it Key words and phrases.\enspace}}\keywords@
 \DNii@##1\endkeywords{\def\thekeywords@{\def\usualspace{{\it\enspace}}%
 \eightpoint\keywords@\ignorespaces##1\unskip.}}%
 \FN@\next@}
\def\thesubjclass@{}
\def\subjclass{\nofrills@{{\rm2010 {\it Mathematics Subject
   Classification\/}.\enspace}}\subjclass@
 \DNii@##1\endsubjclass{\def\thesubjclass@{\def\usualspace
  {{\rm\enspace}}\eightpoint\subjclass@\ignorespaces##1\unskip.}}%
 \FN@\next@}
\newbox\abstractbox@
\def\abstract{\nofrills@{{\smc Abstract.\enspace}}\abstract@
 \DNii@{\setbox\abstractbox@\vbox\bgroup\noindent$$\vbox\bgroup
  \def\envir@{abstract}\advance\hsize-2\indenti
  \usualspace@{{\enspace}}\eightpoint \noindent\abstract@\ignorespaces}%
 \FN@\next@}
\def\endabstract{\par\unskip\egroup$$\egroup}
\def\widestnumber#1#2{\begingroup\let\head\null\let\subhead\empty
   \let\subsubhead\subhead
   \ifx#1\head\global\setbox\tocheadbox@\hbox{#2.\enspace}%
   \else\ifx#1\subhead\global\setbox\tocsubheadbox@\hbox{#2.\enspace}%
   \else\ifx#1\key\bgroup\let\endrefitem@\egroup
        \key#2\endrefitem@\global\refindentwd\wd\keybox@
   \else\ifx#1\no\bgroup\let\endrefitem@\egroup
        \no#2\endrefitem@\global\refindentwd\wd\nobox@
   \else\ifx#1\page\global\setbox\pagesbox@\hbox{\quad\bf#2}%
   \else\ifx#1\item\setboxz@h{#2}\global\rosteritemwd\wdz@
        \global\advance\rosteritemwd by.5\parindent
   \else\message{\string\widestnumber is not defined for this option
   (\string#1)}%
\fi\fi\fi\fi\fi\fi\endgroup}
\newif\ifmonograph@
\def\Monograph{\monograph@true \let\headmark\rightheadtext
  \let\varindent@\indent \def\headfont@{\bf}\def\proclaimheadfont@{\smc}%
  \def\demofont@{\smc}}
%Bei Proclaim,...: Einrcken.
\let\varindent@\indent

\newbox\tocheadbox@    \newbox\tocsubheadbox@
\newbox\tocbox@
\def\toc{\toc@{Contents}}
\def\newtocdefs{%
   \def \title##1\endtitle
       {\penaltyandskip@\z@\smallskipamount
        \hangindent\wd\tocheadbox@\noindent{\bf##1}}%
   \def \chapter##1{%
        Chapter \uppercase\expandafter{\romannumeral##1.\unskip}\enspace}%
   \def \specialhead##1\endspecialhead
       {\par\hangindent\wd\tocheadbox@ \noindent##1\par}%
   \def \head##1 ##2\endhead
       {\par\hangindent\wd\tocheadbox@ \noindent
        \if\notempty{##1}\hbox to\wd\tocheadbox@{\hfil##1\enspace}\fi
        ##2\par}%
   \def \subhead##1 ##2\endsubhead
       {\par\vskip-\parskip {\normalbaselines
        \advance\leftskip\wd\tocheadbox@
        \hangindent\wd\tocsubheadbox@ \noindent
        \if\notempty{##1}\hbox to\wd\tocsubheadbox@{##1\unskip\hfil}\fi
         ##2\par}}%
   \def \subsubhead##1 ##2\endsubsubhead
       {\par\vskip-\parskip {\normalbaselines
        \advance\leftskip\wd\tocheadbox@
        \hangindent\wd\tocsubheadbox@ \noindent
        \if\notempty{##1}\hbox to\wd\tocsubheadbox@{##1\unskip\hfil}\fi
        ##2\par}}}
\def\toc@#1{\relaxnext@
   \def\page##1%
       {\unskip\penalty0\null\hfil
        \rlap{\hbox to\wd\pagesbox@{\quad\hfil##1}}\hfilneg\penalty\@M}%
 \DN@{\ifx\next\nofrills\DN@\nofrills{\nextii@}%
      \else\DN@{\nextii@{{#1}}}\fi
      \next@}%
 \DNii@##1{%
\ifmonograph@\bgroup\else\setbox\tocbox@\vbox\bgroup
   \centerline{\headfont@\ignorespaces##1\unskip}\nobreak
   \vskip\belowheadskip \fi
   \setbox\tocheadbox@\hbox{0.\enspace}%
   \setbox\tocsubheadbox@\hbox{0.0.\enspace}%
   \leftskip\indenti \rightskip\leftskip
   \setbox\pagesbox@\hbox{\bf\quad000}\advance\rightskip\wd\pagesbox@
   \newtocdefs
 }%
 \FN@\next@}
\def\endtoc{\par\egroup}
\let\pretitle\relax
\let\preauthor\relax
\let\preaffil\relax
\let\predate\relax
\let\preabstract\relax
\let\prepaper\relax
\def\dedicatory #1\enddedicatory{\def\preabstract{{\medskip
  \eightpoint\it \raggedcenter@#1\endgraf}}}
\def\thetranslator@{}
\def\translator#1\endtranslator{\def\thetranslator@{\nobreak\medskip
 \line{\eightpoint\hfil Translated by \uppercase{#1}\qquad\qquad}\nobreak}}
\outer\def\endtopmatter{\runaway@{abstract}%
 \edef\next{\the\leftheadtoks}\ifx\next\empty
  \expandafter\leftheadtext\expandafter{\the\rightheadtoks}\fi
 \ifmonograph@\else
   \ifx\thesubjclass@\empty\else \makefootnote@{}{\thesubjclass@}\fi
   \ifx\thekeywords@\empty\else \makefootnote@{}{\thekeywords@}\fi
   \ifx\thethanks@\empty\else \makefootnote@{}{\thethanks@}\fi
 \fi
  \pretitle
  \ifmonograph@ \topskip7pc \else \topskip4pc \fi
  \box\titlebox@
  \topskip10pt% reset to normal value
  \preauthor
  \ifvoid\authorbox@\else \vskip2.5pc plus1pc \unvbox\authorbox@\fi
  \preaffil
  \ifvoid\affilbox@\else \vskip1pc plus.5pc \unvbox\affilbox@\fi
  \predate
  \ifx\thedate@\empty\else \vskip1pc plus.5pc \line{\hfil\thedate@\hfil}\fi
  \preabstract
  \ifvoid\abstractbox@\else \vskip1.5pc plus.5pc \unvbox\abstractbox@ \fi
  \ifvoid\tocbox@\else\vskip1.5pc plus.5pc \unvbox\tocbox@\fi
  \prepaper
  \vskip2pc plus1pc
}
\def\document{\let\fontlist@\relax\let\alloclist@\relax
  \tenpoint}

%Modifizierte Head-Skips
\newskip\aboveheadskip       \aboveheadskip1.8\bigskipamount
\newdimen\belowheadskip      \belowheadskip1.8\medskipamount

\def\headfont@{\smc}
\def\penaltyandskip@#1#2{\relax\ifdim\lastskip<#2\relax\removelastskip
      \ifnum#1=\z@\else\penalty@#1\relax\fi\vskip#2%
  \else\ifnum#1=\z@\else\penalty@#1\relax\fi\fi}
\def\nobreak{\penalty\@M
  \ifvmode\def\penalty@{\let\penalty@\penalty\count@@@}%
  \everypar{\let\penalty@\penalty\everypar{}}\fi}
\let\penalty@\penalty
\def\heading#1\endheading{\head#1\endhead}

\def\specialheadfont@{\bf}
\outer\def\specialhead{\par\penaltyandskip@{-200}\aboveheadskip
  \begingroup\interlinepenalty\@M\rightskip\z@ plus\hsize \let\\\linebreak
  \specialheadfont@\noindent\ignorespaces}
\def\endspecialhead{\par\endgroup\nobreak\vskip\belowheadskip}
%\outer\def\head#1\endhead{\par\penaltyandskip@{-200}\aboveheadskip
% {\headfont@\raggedcenter@\interlinepenalty\@M
% \ignorespaces#1\endgraf}\nobreak
% \vskip\belowheadskip
% \headmark{#1}}
\let\headmark\eat@
\newskip\subheadskip       \subheadskip\medskipamount
\def\subheadfont@{\bf}
\outer\def\subhead{\nofrills@{.\enspace}\subhead@
 \DNii@##1\endsubhead{\par\penaltyandskip@{-100}\subheadskip
  \varindent@{\usualspace@{{\subheadfont@\enspace}}%
 \subheadfont@\ignorespaces##1\unskip\subhead@}\ignorespaces}%
 \FN@\next@}
\outer\def\subsubhead{\nofrills@{.\enspace}\subsubhead@
 \DNii@##1\endsubsubhead{\par\penaltyandskip@{-50}\medskipamount
      {\usualspace@{{\it\enspace}}%
  \it\ignorespaces##1\unskip\subsubhead@}\ignorespaces}%
 \FN@\next@}
\def\proclaimheadfont@{\bf}
\outer\def\proclaim{\runaway@{proclaim}\def\envir@{proclaim}%
  \nofrills@{.\enspace}\proclaim@
 \DNii@##1{\penaltyandskip@{-100}\medskipamount\varindent@
   \usualspace@{{\proclaimheadfont@\enspace}}\proclaimheadfont@
   \ignorespaces##1\unskip\proclaim@
  \sl\ignorespaces}% 
 \FN@\next@}
\outer\def\endproclaim{\let\envir@\relax\par\rm
  \penaltyandskip@{55}\medskipamount}
\def\demoheadfont@{\it}
\def\demo{\runaway@{proclaim}\nofrills@{.\enspace}\demo@
     \DNii@##1{\par\penaltyandskip@\z@\medskipamount
  {\usualspace@{{\demoheadfont@\enspace}}%
  \varindent@\demoheadfont@\ignorespaces##1\unskip\demo@}\rm
  \ignorespaces}\FN@\next@}
\def\enddemo{\par\medskip}
\def\qed{\ifhmode\unskip\nobreak\fi\quad\ifmmode\square\else$\m@th\square$\fi}
\let\remark\demo
%\endremark=\enddemo
\let\endremark\enddemo
%Headfont fr Definition wie bei Beweisen.
\def\definition{\runaway@{proclaim}%
  \nofrills@{.\demoheadfont@\enspace}\definition@
        \DNii@##1{\penaltyandskip@{-100}\medskipamount
        {\usualspace@{{\demoheadfont@\enspace}}%
        \varindent@\demoheadfont@\ignorespaces##1\unskip\definition@}%
        \rm \ignorespaces}\FN@\next@}

%\let\example\definition
%\let\endexample\enddefinition
%Modifikation:

\newdimen\rosteritemwd
\newcount\rostercount@
\newif\iffirstitem@
\let\plainitem@\item
\newtoks\everypartoks@
\def\par@{\everypartoks@\expandafter{\the\everypar}\everypar{}}
\def\roster{\edef\leftskip@{\leftskip\the\leftskip}%
 \relaxnext@
 \rostercount@\z@  
 \def\item{\FN@\rosteritem@}% 
 \DN@{\ifx\next\runinitem\let\next@\nextii@\else
  \let\next@\nextiii@\fi\next@}%
 \DNii@\runinitem  
  {\unskip  
   \DN@{\ifx\next[\let\next@\nextii@\else
    \ifx\next"\let\next@\nextiii@\else\let\next@\nextiv@\fi\fi\next@}%
   \DNii@[####1]{\rostercount@####1\relax
    \enspace{\rm(\number\rostercount@)}~\ignorespaces}%
   \def\nextiii@"####1"{\enspace{\rm####1}~\ignorespaces}%
   \def\nextiv@{\enspace{\rm(1)}\rostercount@\@ne~}%
   \par@\firstitem@false  
   \FN@\next@}% 
 \def\nextiii@{\par\par@  
  \penalty\@m\smallskip\vskip-\parskip
  \firstitem@true}%  
 \FN@\next@}
\def\rosteritem@{\iffirstitem@\firstitem@false\else\par\vskip-\parskip\fi
 \leftskip3\parindent\noindent  
 \DNii@[##1]{\rostercount@##1\relax
  \llap{\hbox to2.5\parindent{\hss\rm(\number\rostercount@)}%
   \hskip.5\parindent}\ignorespaces}%
 \def\nextiii@"##1"{%
  \llap{\hbox to2.5\parindent{\hss\rm##1}\hskip.5\parindent}\ignorespaces}%
 \def\nextiv@{\advance\rostercount@\@ne
  \llap{\hbox to2.5\parindent{\hss\rm(\number\rostercount@)}%
   \hskip.5\parindent}}%
 \ifx\next[\let\next@\nextii@\else\ifx\next"\let\next@\nextiii@\else
  \let\next@\nextiv@\fi\fi\next@}

\newif\ifnextRunin@
\def\endroster{\relaxnext@
 \par\leftskip@  
 \penalty-50 \vskip-\parskip\smallskip  
 \DN@{\ifx\next\Runinitem\let\next@\relax
  \else\nextRunin@false\let\item\plainitem@  
   \ifx\next\par 
    \DN@\par{\everypar\expandafter{\the\everypartoks@}}%
   \else  
    \DN@{\noindent\everypar\expandafter{\the\everypartoks@}}%
  \fi\fi\next@}%
 \FN@\next@}
\newcount\rosterhangafter@
\def\Runinitem#1\roster\runinitem{\relaxnext@
 \rostercount@\z@ 
 \def\item{\FN@\rosteritem@}%  
 \def\runinitem@{#1}%  
 \DN@{\ifx\next[\let\next\nextii@\else\ifx\next"\let\next\nextiii@
  \else\let\next\nextiv@\fi\fi\next}%
 \DNii@[##1]{\rostercount@##1\relax
  \def\item@{{\rm(\number\rostercount@)}}\nextv@}%
 \def\nextiii@"##1"{\def\item@{{\rm##1}}\nextv@}%
 \def\nextiv@{\advance\rostercount@\@ne
  \def\item@{{\rm(\number\rostercount@)}}\nextv@}%
 \def\nextv@{\setbox\z@\vbox  
  {\ifnextRunin@\noindent\fi  
  \runinitem@\unskip\enspace\item@~\par  
  \global\rosterhangafter@\prevgraf}% 
  \firstitem@false  
  \ifnextRunin@\else\par\fi  
  \hangafter\rosterhangafter@\hangindent3\parindent
  \ifnextRunin@\noindent\fi  
  \runinitem@\unskip\enspace 
  \item@~\ifnextRunin@\else\par@\fi  
  \nextRunin@true\ignorespaces}%  
 \FN@\next@}
\def\footmarkform@#1{$\m@th^{#1}$}
\let\thefootnotemark\footmarkform@
\def\makefootnote@#1#2{\insert\footins
 {\interlinepenalty\interfootnotelinepenalty
 \eightpoint\splittopskip\ht\strutbox\splitmaxdepth\dp\strutbox
 \floatingpenalty\@MM\leftskip\z@\rightskip\z@\spaceskip\z@\xspaceskip\z@
 \leavevmode{#1}\footstrut\ignorespaces#2\unskip\lower\dp\strutbox
 \vbox to\dp\strutbox{}}}
\newcount\footmarkcount@
\footmarkcount@\z@
\def\footnotemark{\let\@sf\empty\relaxnext@
 \ifhmode\edef\@sf{\spacefactor\the\spacefactor}\/\fi
 \DN@{\ifx[\next\let\next@\nextii@\else
  \ifx"\next\let\next@\nextiii@\else
  \let\next@\nextiv@\fi\fi\next@}%
 \DNii@[##1]{\footmarkform@{##1}\@sf}%
 \def\nextiii@"##1"{{##1}\@sf}%
 \def\nextiv@{\iffirstchoice@\global\advance\footmarkcount@\@ne\fi
  \footmarkform@{\number\footmarkcount@}\@sf}%
 \FN@\next@}
\def\footnotetext{\relaxnext@
 \DN@{\ifx[\next\let\next@\nextii@\else
  \ifx"\next\let\next@\nextiii@\else
  \let\next@\nextiv@\fi\fi\next@}%
 \DNii@[##1]##2{\makefootnote@{\footmarkform@{##1}}{##2}}%
 \def\nextiii@"##1"##2{\makefootnote@{##1}{##2}}%
 \def\nextiv@##1{\makefootnote@{\footmarkform@{\number\footmarkcount@}}{##1}}%
 \FN@\next@}
\def\footnote{\let\@sf\empty\relaxnext@
 \ifhmode\edef\@sf{\spacefactor\the\spacefactor}\/\fi
 \DN@{\ifx[\next\let\next@\nextii@\else
  \ifx"\next\let\next@\nextiii@\else
  \let\next@\nextiv@\fi\fi\next@}%
 \DNii@[##1]##2{\footnotemark[##1]\footnotetext[##1]{##2}}%
 \def\nextiii@"##1"##2{\footnotemark"##1"\footnotetext"##1"{##2}}%
 \def\nextiv@##1{\footnotemark\footnotetext{##1}}%
 \FN@\next@}
\def\adjustfootnotemark#1{\advance\footmarkcount@#1\relax}
\def\footnoterule{\kern-3\p@
  \hrule width 5pc\kern 2.6\p@} 
\def\captionfont@{\smc}
\def\topcaption#1#2\endcaption{%
  {\dimen@\hsize \advance\dimen@-\captionwidth@
   \rm\raggedcenter@ \advance\leftskip.5\dimen@ \rightskip\leftskip
  {\captionfont@#1}%
  \if\notempty{#2}.\enspace\ignorespaces#2\fi
  \endgraf}\nobreak\bigskip}
\def\botcaption#1#2\endcaption{%
  \nobreak\bigskip
  \setboxz@h{\captionfont@#1\if\notempty{#2}.\enspace\rm#2\fi}%
  {\dimen@\hsize \advance\dimen@-\captionwidth@
   \leftskip.5\dimen@ \rightskip\leftskip
   \noindent \ifdim\wdz@>\captionwidth@ 
   \else\hfil\fi 
  {\captionfont@#1}\if\notempty{#2}.\enspace\rm#2\fi\endgraf}}
\def\@ins{\par\begingroup\def\vspace##1{\vskip##1\relax}%
  \def\captionwidth##1{\captionwidth@##1\relax}%
  \setbox\z@\vbox\bgroup} % start a \vbox
\def\block{\RIfMIfI@\nondmatherr@\block\fi
       \else\ifvmode\vskip\abovedisplayskip\noindent\fi
        $$\def\endblock{\par\egroup$$}\fi
  \vbox\bgroup\advance\hsize-2\indenti\noindent}
\def\endblock{\par\egroup}
\def\cite#1{{\rm[{\citefont@\m@th#1}]}}
\def\citefont@{\rm}
\def\refsfont@{\eightpoint}
\outer\def\Refs{\runaway@{proclaim}%
 \relaxnext@ \DN@{\ifx\next\nofrills\DN@\nofrills{\nextii@}\else
  \DN@{\nextii@{References}}\fi\next@}%
 \DNii@##1{\penaltyandskip@{-200}\aboveheadskip
  \line{\hfil\headfont@\ignorespaces##1\unskip\hfil}\nobreak
  \vskip\belowheadskip
  \begingroup\refsfont@\sfcode`.=\@m}%
 \FN@\next@}
\def\endRefs{\par\endgroup}
\newbox\nobox@            \newbox\keybox@           \newbox\bybox@
\newbox\paperbox@         \newbox\paperinfobox@     \newbox\jourbox@
\newbox\volbox@           \newbox\issuebox@         \newbox\yrbox@
\newbox\pagesbox@         \newbox\bookbox@          \newbox\bookinfobox@
\newbox\publbox@          \newbox\publaddrbox@      \newbox\finalinfobox@
\newbox\edsbox@           \newbox\langbox@
\newif\iffirstref@        \newif\iflastref@
\newif\ifprevjour@        \newif\ifbook@            \newif\ifprevinbook@
\newif\ifquotes@          \newif\ifbookquotes@      \newif\ifpaperquotes@
\newdimen\bysamerulewd@
\setboxz@h{\refsfont@\kern3em}
\bysamerulewd@\wdz@
\newdimen\refindentwd
\setboxz@h{\refsfont@ 00. }
\refindentwd\wdz@
\outer\def\ref{\begingroup \noindent\hangindent\refindentwd
 \firstref@true \def\nofrills{\def\refkern@{\kern3sp}}%
 \ref@}
\def\ref@{\book@false \bgroup\let\endrefitem@\egroup \ignorespaces}
\def\moreref{\endrefitem@\endref@\firstref@false\ref@}%
\def\transl{\endrefitem@\endref@\firstref@false
  \book@false
  \prepunct@
  \setboxz@h\bgroup \aftergroup\unhbox\aftergroup\z@
    \def\endrefitem@{\unskip\refkern@\egroup}\ignorespaces}%
\def\emptyifempty@{\dimen@\wd\currbox@
  \advance\dimen@-\wd\z@ \advance\dimen@-.1\p@
  \ifdim\dimen@<\z@ \setbox\currbox@\copy\voidb@x \fi}
\let\refkern@\relax
\def\endrefitem@{\unskip\refkern@\egroup
  \setboxz@h{\refkern@}\emptyifempty@}\ignorespaces
\def\refdef@#1#2#3{\edef\next@{\noexpand\endrefitem@
  \let\noexpand\currbox@\csname\expandafter\eat@\string#1box@\endcsname
    \noexpand\setbox\noexpand\currbox@\hbox\bgroup}%
  \toks@\expandafter{\next@}%
  \if\notempty{#2#3}\toks@\expandafter{\the\toks@
  \def\endrefitem@{\unskip#3\refkern@\egroup
  \setboxz@h{#2#3\refkern@}\emptyifempty@}#2}\fi
  \toks@\expandafter{\the\toks@\ignorespaces}%
  \edef#1{\the\toks@}}
\refdef@\no{}{. }
\refdef@\key{[\m@th}{] }
\refdef@\by{}{}
\def\bysame{\by\hbox to\bysamerulewd@{\hrulefill}\thinspace
   \kern0sp}
\def\manyby{\message{\string\manyby is no longer necessary; \string\by
  can be used instead, starting with version 2.0 of \styname.STY}\by}
\refdef@\paper{\ifpaperquotes@``\fi\it}{}
\refdef@\paperinfo{}{}
\def\jour{\endrefitem@\let\currbox@\jourbox@
  \setbox\currbox@\hbox\bgroup
  \def\endrefitem@{\unskip\refkern@\egroup
    \setboxz@h{\refkern@}\emptyifempty@
    \ifvoid\jourbox@\else\prevjour@true\fi}%
\ignorespaces}
\refdef@\vol{\ifbook@\else\bf\fi}{}
\refdef@\issue{no. }{}
\refdef@\yr{}{}
\refdef@\pages{}{}
\def\page{\endrefitem@\def\pp@{\def\pp@{pp.~}p.~}\let\currbox@\pagesbox@
  \setbox\currbox@\hbox\bgroup\ignorespaces}
\def\pp@{pp.~}
\def\book{\endrefitem@ \let\currbox@\bookbox@
 \setbox\currbox@\hbox\bgroup\def\endrefitem@{\unskip\refkern@\egroup
  \setboxz@h{\ifbookquotes@``\fi}\emptyifempty@
  \ifvoid\bookbox@\else\book@true\fi}%
  \ifbookquotes@``\fi\it\ignorespaces}
\def\inbook{\endrefitem@
  \let\currbox@\bookbox@\setbox\currbox@\hbox\bgroup
  \def\endrefitem@{\unskip\refkern@\egroup
  \setboxz@h{\ifbookquotes@``\fi}\emptyifempty@
  \ifvoid\bookbox@\else\book@true\previnbook@true\fi}%
  \ifbookquotes@``\fi\ignorespaces}
\refdef@\eds{(}{, eds.)}
\def\ed{\endrefitem@\let\currbox@\edsbox@
 \setbox\currbox@\hbox\bgroup
 \def\endrefitem@{\unskip, ed.)\refkern@\egroup
  \setboxz@h{(, ed.)}\emptyifempty@}(\ignorespaces}
\refdef@\bookinfo{}{}
\refdef@\publ{}{}
\refdef@\publaddr{}{}
\refdef@\finalinfo{}{}
\refdef@\lang{(}{)}

\let\refdef@\relax 
\def\ppunbox@#1{\ifvoid#1\else\prepunct@\unhbox#1\fi}
\def\nocomma@#1{\ifvoid#1\else\changepunct@3\prepunct@\unhbox#1\fi}
\def\changepunct@#1{\ifnum\lastkern<3 \unkern\kern#1sp\fi}
\def\prepunct@{\count@\lastkern\unkern
  \ifnum\lastpenalty=0
    \let\penalty@\relax
  \else
    \edef\penalty@{\penalty\the\lastpenalty\relax}%
  \fi
  \unpenalty
  \let\refspace@\ \ifcase\count@,% usual case, do a comma
\or;\or.\or % do nothing; this case is from nofrills.
  \or\let\refspace@\relax
  \else,\fi
  \ifquotes@''\quotes@false\fi \penalty@ \refspace@
}
\def\transferpenalty@#1{\dimen@\lastkern\unkern
  \ifnum\lastpenalty=0\unpenalty\let\penalty@\relax
  \else\edef\penalty@{\penalty\the\lastpenalty\relax}\unpenalty\fi
  #1\penalty@\kern\dimen@}
\def\endref{\endrefitem@\lastref@true\endref@
  \par\endgroup \prevjour@false \previnbook@false }
\def\endref@{%
\iffirstref@
  \ifvoid\nobox@\ifvoid\keybox@\indent\fi
  \else\hbox to\refindentwd{\hss\unhbox\nobox@}\fi
  \ifvoid\keybox@
  \else\ifdim\wd\keybox@>\refindentwd
         \box\keybox@
       \else\hbox to\refindentwd{\unhbox\keybox@\hfil}\fi\fi
  \kern4sp\ppunbox@\bybox@
\fi 
  \ifvoid\paperbox@
  \else\prepunct@\unhbox\paperbox@
    \ifpaperquotes@\quotes@true\fi\fi
  \ppunbox@\paperinfobox@
  \ifvoid\jourbox@
    \ifprevjour@ \nocomma@\volbox@
      \nocomma@\issuebox@
      \ifvoid\yrbox@\else\changepunct@3\prepunct@(\unhbox\yrbox@
        \transferpenalty@)\fi
      \ppunbox@\pagesbox@
    \fi 
  \else \prepunct@\unhbox\jourbox@
    \nocomma@\volbox@
    \nocomma@\issuebox@
    \ifvoid\yrbox@\else\changepunct@3\prepunct@(\unhbox\yrbox@
      \transferpenalty@)\fi
    \ppunbox@\pagesbox@
  \fi 
  \ifbook@\prepunct@\unhbox\bookbox@ \ifbookquotes@\quotes@true\fi \fi
  \nocomma@\edsbox@
  \ppunbox@\bookinfobox@
  \ifbook@\ifvoid\volbox@\else\prepunct@ vol.~\unhbox\volbox@
  \fi\fi
  \ppunbox@\publbox@ \ppunbox@\publaddrbox@
  \ifbook@ \ppunbox@\yrbox@
    \ifvoid\pagesbox@
    \else\prepunct@\pp@\unhbox\pagesbox@\fi
  \else
    \ifprevinbook@ \ppunbox@\yrbox@
      \ifvoid\pagesbox@\else\prepunct@\pp@\unhbox\pagesbox@\fi
    \fi \fi
  \ppunbox@\finalinfobox@
  \iflastref@
    \ifvoid\langbox@.\ifquotes@''\fi
    \else\changepunct@2\prepunct@\unhbox\langbox@\fi
  \else
    \ifvoid\langbox@\changepunct@1%
    \else\changepunct@3\prepunct@\unhbox\langbox@
      \changepunct@1\fi
  \fi
}
\outer\def\enddocument{%
 \runaway@{proclaim}%
\ifmonograph@ % do nothing
\else
 \nobreak
 \thetranslator@
 \count@\z@ \loop\ifnum\count@<\addresscount@\advance\count@\@ne
 \csname address\number\count@\endcsname
 \csname email\number\count@\endcsname
 \repeat
\fi
 \vfill\supereject\end}

%Modifizierte Fonts fr Proclaim,...
\def\headfont@{\headfonts}
\def\proclaimheadfont@{\bf}
\def\specialheadfont@{\bf}
\def\subheadfont@{\bf}
\def\demoheadfont@{\smc}

%Kontrollsequenzen fr Inhaltsverzeichnis und Index:
\newif\ifThisToToc \ThisToTocfalse
\newif\iftocloaded \tocloadedfalse

\def\C@L{\noexpand\Cal}\def\B@B{\noexpand\Bbb}\def\fR@K{\noexpand\frak}
\def\S@{\noexpand\S}\def\P@P{\noexpand\"}
\def\xpar{\\}

\def\writetoc#1{\iftocloaded\ifThisToToc\begingroup\def\totoc{}
  \def\Cal{\noexpand\C@L}\def\Bbb{\noexpand\B@B}
  \def\frak{\noexpand\fR@K}\def\goth{\frak}\def\S{\noexpand\S@}
  \def\"{\noexpand\P@P}
  \def\xpar{\par\penalty100000 }\def\idx##1{##1}\def\\{\xpar}
  \edef\next@{\write\toc{\noindent#1\leaderfill\noexpand\folio\par}}%
  \next@\endgroup\global\ThisToTocfalse\fi\fi}
\def\leaderfill{\leaders\hbox to 1em{\hss.\hss}\hfill}

\newif\ifindexloaded \indexloadedfalse
\def\idx#1{\ifindexloaded\begingroup\def\ign{}\def\it{}\def\/{}%
 \def\smc{}\def\bf{}\def\tt{}%
 \def\Cal{\noexpand\C@L}\def\Bbb{\noexpand\B@B}%
 \def\frak{\noexpand\fR@K}\def\goth{\frak}\def\S{\noexpand\S@}%
  \def\"{\noexpand\P@P}%
 {\edef\next@{\write\index{#1, \noexpand\folio}}\next@}%
 \endgroup\fi{#1}}
\def\ign#1{}

\def\totoc{\global\ThisToToctrue}

\outer\def\head#1\endhead{\par\penaltyandskip@{-200}\aboveheadskip
 {\headfont@\raggedcenter@\interlinepenalty\@M
 \ignorespaces#1\endgraf}\nobreak
 \vskip\belowheadskip
 \headmark{#1}\writetoc{#1}}

\outer\def\chaphead#1\endchaphead{\par\penaltyandskip@{-200}\aboveheadskip
 {\chapheadfonts\raggedcenter@\interlinepenalty\@M
 \ignorespaces#1\endgraf}\nobreak
 \vskip3\belowheadskip
 \headmark{#1}\writetoc{#1}}

\def\folio{{\foliofont@\ifnum\pageno<\z@ \romannumeral-\pageno
 \else\number\pageno \fi}}
\newtoks\leftheadtoks
\newtoks\rightheadtoks

%Uppercase ist abgestellt:
\def\leftheadtext{\nofrills@{\relax}\lht@
  \DNii@##1{\leftheadtoks\expandafter{\lht@{##1}}%
    \mark{\the\leftheadtoks\noexpand\else\the\rightheadtoks}
    \ifsyntax@\setboxz@h{\def\\{\unskip\space\ignorespaces}%
        \headlinefont@##1}\fi}%
  \FN@\next@}
%Uppercase ist abgestellt:
\def\rightheadtext{\nofrills@{\relax}\rht@
  \DNii@##1{\rightheadtoks\expandafter{\rht@{##1}}%
    \mark{\the\leftheadtoks\noexpand\else\the\rightheadtoks}%
    \ifsyntax@\setboxz@h{\def\\{\unskip\space\ignorespaces}%
        \headlinefont@##1}\fi}%
  \FN@\next@}
%\headline={\def\chapter#1{\chapterno@. }%
%  \def\\{\unskip\space\ignorespaces}\headlinefont@
%  \ifodd\pageno \rightheadline \else \leftheadline\fi}
\def\NoRunningHeads{\global\runheads@false\global\let\headmark\eat@}

\newif\iffirstpage@     \firstpage@true
\newif\ifrunheads@      \runheads@true

%Ergnzungen zu Runningheads und Pagenumbers:
\newdimen\fullhsize \fullhsize=\hsize
\newdimen\fullvsize \fullvsize=\vsize
\def\fullline{\hbox to\fullhsize}

\def\pagenumbers{\gdef\folio{\folio@}}

\let\norunningheads\NoRunningHeads
\def\userunningheads{\global\runheads@true}
%Default: Seitennumerierung unten.
\norunningheads

\headline={\def\chapter#1{\chapterno@. }%
  \def\\{\unskip\space\ignorespaces}\ifrunheads@\headlinefont@
    \ifodd\pageno\rightheadline \else\leftheadline\fi
   \else\hfil\fi\ifNoRunHeadline\global\NoRunHeadlinefalse\fi}
\let\folio@\folio
\def\foliofont@{\foliofont}
\def\foliofont{\eightrm}
\def\headlinefont@{\headlinefont}
\def\headlinefont{\eightpoint\smc}
\def\leftheadline{\rlap{\folio}\hfill
   \ifNoRunHeadline\else\iftrue\topmark\fi\fi \hfill}
\def\rightheadline{\hfill\ifNoRunHeadline
   \else \expandafter\iffalse\botmark\fi\fi
  \hfill \llap{\folio}}
\footline={{\eightpoint\bottremark}%
   \ifrunheads@\else\hfil{\let\foliofont\tenrm\folio}\fi\hfil}
\def\bottremark{}
 
%Definition von \norunninghead:
\newif\ifNoRunHeadline      
\def\norunninghead{\global\NoRunHeadlinetrue}
\norunninghead

\output={\output@}
%\def\output@{\shipout\vbox{%
% \iffirstpage@ \global\firstpage@false
%  \pagebody \logo@ \makefootline%
% \else \ifrunheads@ \makeheadline \pagebody
%       \else \pagebody \makefootline \fi
% \fi}%
% \advancepageno \ifnum\outputpenalty>-\@MM\else\dosupereject\fi}
%
%Modifizierter Output + Index-Output:
\newif\ifoffset\offsetfalse
\output={\output@}
\def\output@{%
 \ifoffset 
  \ifodd\count0\advance\hoffset by0.5truecm
   \else\advance\hoffset by-0.5truecm\fi\fi
 \shipout\vbox{%
  \makeheadline \pagebody \makefootline }%
 \advancepageno \ifnum\outputpenalty>-\@MM\else\dosupereject\fi}

\def\indexoutput#1{%
  \ifoffset 
   \ifodd\count0\advance\hoffset by0.5truecm
    \else\advance\hoffset by-0.5truecm\fi\fi
  \shipout\vbox{\makeheadline
  \vbox to\fullvsize{\boxmaxdepth\maxdepth%
  \ifvoid\topins\else\unvbox\topins\fi% 
  #1 %
  \ifvoid\footins\else % footnote info is present
    \vskip\skip\footins
    \footnoterule
    \unvbox\footins\fi
  \ifr@ggedbottom \kern-\dimen@ \vfil \fi}%
  \baselineskip2pc
  \makefootline}%
 \global\advance\pageno\@ne
 \ifnum\outputpenalty>-\@MM\else\dosupereject\fi}
 
 \newbox\partialpage \newdimen\halfsize \halfsize=0.5\fullhsize
 \advance\halfsize by-0.5em

 \def\begindoublecolumns{\output={\indexoutput{\unvbox255}}%
   \begingroup \def\line{\fullline}
   \output={\global\setbox\partialpage=\vbox{\unvbox255\bigskip}}\eject
   \output={\doublecolumnout}\hsize=\halfsize \vsize=2\fullvsize}
 \def\enddoublecolumns{\output={\balancecolumns}\eject
  \endgroup \pagegoal=\fullvsize%
  \output={\output@}}
\def\doublecolumnout{\splittopskip=\topskip \splitmaxdepth=\maxdepth
  \dimen@=\fullvsize \advance\dimen@ by-\ht\partialpage
  \setbox0=\vsplit255 to \dimen@ \setbox2=\vsplit255 to \dimen@
  \indexoutput{\pagesofar} \unvbox255 \penalty\outputpenalty}
\def\pagesofar{\unvbox\partialpage
  \wd0=\hsize \wd2=\hsize \hbox to\fullhsize{\box0\hfil\box2}}
\def\balancecolumns{\setbox0=\vbox{\unvbox255} \dimen@=\ht0
  \advance\dimen@ by\topskip \advance\dimen@ by-\baselineskip
  \divide\dimen@ by2 \splittopskip=\topskip
  {\vbadness=10000 \loop \global\setbox3=\copy0
    \global\setbox1=\vsplit3 to\dimen@
    \ifdim\ht3>\dimen@ \global\advance\dimen@ by1pt \repeat}
  \setbox0=\vbox to\dimen@{\unvbox1} \setbox2=\vbox to\dimen@{\unvbox3}
  \pagesofar}

\tenpoint
\catcode`\@=\active

\def\smallheadings{\let\chapheadfonts\tenpoint\let\headfonts\tenpoint}

\tenpoint
\catcode`\@=\active

\def\LL{\leavevmode\setbox0=\hbox{L}\hbox to\wd0{\hss\char'40L}}
\def\al{\alpha}
\def\be{\beta}
\def\ga{\gamma}
\def\de{\delta}
\def\ep{\varepsilon}

\def\th{\theta}

\def\la{\lambda}

\def\si{\sigma}

            %used for crossreferencing, Tex should ignore.
             %used for refencing (section-numbers)
          %used for new-section numbers

\def\Z{{\Bbb Z}}

\def\today{\ifcase\month\or
 January\or February\or March\or April\or May\or June\or
 July\or August\or September\or October\or November\or December\fi
 \space\number\day, \number\year}
 %zum Nummerieren

\def\({\left(}
\def\){\right)}
\def\[{\left[}
\def\]{\right]}

\def\sgn{\operatorname{sgn}}

\def\3{\ss}
\catcode`\@=11
\def\dddot#1{\vbox{\ialign{##\crcr
      .\hskip-.5pt.\hskip-.5pt.\crcr\noalign{\kern1.5\p@\nointerlineskip}
      $\hfil\displaystyle{#1}\hfil$\crcr}}}

\newif\iftab@\tab@false
\newif\ifvtab@\vtab@false
\def\tab{\bgroup\tab@true\vtab@false\vst@bfalse\Strich@false%
   \def\\{\global\hline@@false%
     \ifhline@\global\hline@false\global\hline@@true\fi\cr}
   \edef\l@{\the\leftskip}\ialign\bgroup\hskip\l@##\hfil&&##\hfil\cr}
\def\endtab{\cr\egroup\egroup}
\def\vtab{\vtop\bgroup\vst@bfalse\vtab@true\tab@true\Strich@false%
   \bgroup\def\\{\cr}\ialign\bgroup&##\hfil\cr}
\def\endvtab{\cr\egroup\egroup\egroup}
\def\stab{\D@cke0.5pt\null 
 \bgroup\tab@true\vtab@false\vst@bfalse\Strich@true\Let@@\vspace@
 \normalbaselines\offinterlineskip
  \openup\spreadmlines@
 \edef\l@{\the\leftskip}\ialign
 \bgroup\hskip\l@##\hfil&&##\hfil\crcr}
\def\endstab{\crcr\egroup
 \egroup}
\newif\ifvst@b\vst@bfalse
\def\vstab{\D@cke0.5pt\null
 \vtop\bgroup\tab@true\vtab@false\vst@btrue\Strich@true\bgroup\Let@@\vspace@
 \normalbaselines\offinterlineskip
  \openup\spreadmlines@\bgroup}
\def\endvstab{\crcr\egroup\egroup
 \egroup\tab@false\Strich@false}

\newdimen\htstrut@
\htstrut@8.5\p@
\newdimen\htStrut@
\htStrut@12\p@
\newdimen\dpstrut@
\dpstrut@3.5\p@
\newdimen\dpStrut@
\dpStrut@3.5\p@
\def\openup{\afterassignment\@penup\dimen@=}
\def\@penup{\advance\lineskip\dimen@
  \advance\baselineskip\dimen@
  \advance\lineskiplimit\dimen@
  \divide\dimen@ by2
  \advance\htstrut@\dimen@
  \advance\htStrut@\dimen@
  \advance\dpstrut@\dimen@
  \advance\dpStrut@\dimen@}
\def\Let@@{\relax\iffalse{\fi%
    \def\\{\global\hline@@false%
     \ifhline@\global\hline@false\global\hline@@true\fi\cr}%
    \iffalse}\fi}
\def\matrix{\null\,\vcenter\bgroup
 \tab@false\vtab@false\vst@bfalse\Strich@false\Let@@\vspace@
 \normalbaselines\openup\spreadmlines@\ialign
 \bgroup\hfil$\m@th##$\hfil&&\quad\hfil$\m@th##$\hfil\crcr
 \Mathstrut@\crcr\noalign{\kern-\baselineskip}}
\def\endmatrix{\crcr\Mathstrut@\crcr\noalign{\kern-\baselineskip}\egroup
 \egroup\,}
\def\smatrix{\D@cke0.5pt\null\,
 \vcenter\bgroup\tab@false\vtab@false\vst@bfalse\Strich@true\Let@@\vspace@
 \normalbaselines\offinterlineskip
  \openup\spreadmlines@\ialign
 \bgroup\hfil$\m@th##$\hfil&&\quad\hfil$\m@th##$\hfil\crcr}
\def\endsmatrix{\crcr\egroup
 \egroup\,\Strich@false}
\newdimen\D@cke
\def\Dicke#1{\global\D@cke#1}
\newtoks\tabs@\tabs@{&}
\newif\ifStrich@\Strich@false
\newif\iff@rst

\def\Stricherr@{\iftab@\ifvtab@\errmessage{\noexpand\s not allowed
     here. Use \noexpand\vstab!}%
  \else\errmessage{\noexpand\s not allowed here. Use \noexpand\stab!}%
  \fi\else\errmessage{\noexpand\s not allowed
     here. Use \noexpand\smatrix!}\fi}
\def\format{\ifvst@b\else\crcr\fi\egroup\iffalse{\fi\ifnum`}=0 \fi\format@}
\def\format@#1\\{\def\preamble@{#1}%
 \def\Str@chfehlt##1{\ifx##1\s\Stricherr@\fi\ifx##1\\\let\Next\relax%
   \else\let\Next\Str@chfehlt\fi\Next}%
 \def\c{\hfil\noexpand\ifhline@@\hbox{\vrule height\htStrut@%
   depth\dpstrut@ width\z@}\noexpand\fi%
   \ifStrich@\hbox{\vrule height\htstrut@ depth\dpstrut@ width\z@}%
   \fi\iftab@\else$\m@th\fi\the\hashtoks@\iftab@\else$\fi\hfil}%
 \def\r{\hfil\noexpand\ifhline@@\hbox{\vrule height\htStrut@%
   depth\dpstrut@ width\z@}\noexpand\fi%
   \ifStrich@\hbox{\vrule height\htstrut@ depth\dpstrut@ width\z@}%
   \fi\iftab@\else$\m@th\fi\the\hashtoks@\iftab@\else$\fi}%
 \def\l{\noexpand\ifhline@@\hbox{\vrule height\htStrut@%
   depth\dpstrut@ width\z@}\noexpand\fi%
   \ifStrich@\hbox{\vrule height\htstrut@ depth\dpstrut@ width\z@}%
   \fi\iftab@\else$\m@th\fi\the\hashtoks@\iftab@\else$\fi\hfil}%
 \def\s{\ifStrich@\ \the\tabs@\vrule width\D@cke\the\hashtoks@%
          \fi\the\tabs@\ }%
 \def\sa{\ifStrich@\vrule width\D@cke\the\hashtoks@%
            \the\tabs@\ %
            \fi}%
 \def\se{\ifStrich@\ \the\tabs@\vrule width\D@cke\the\hashtoks@\fi}%
 \def\cd{\hfil\noexpand\ifhline@@\hbox{\vrule height\htStrut@%
   depth\dpstrut@ width\z@}\noexpand\fi%
   \ifStrich@\hbox{\vrule height\htstrut@ depth\dpstrut@ width\z@}%
   \fi$\dsize\m@th\the\hashtoks@$\hfil}%
 \def\rd{\hfil\noexpand\ifhline@@\hbox{\vrule height\htStrut@%
   depth\dpstrut@ width\z@}\noexpand\fi%
   \ifStrich@\hbox{\vrule height\htstrut@ depth\dpstrut@ width\z@}%
   \fi$\dsize\m@th\the\hashtoks@$}%
 \def\ld{\noexpand\ifhline@@\hbox{\vrule height\htStrut@%
   depth\dpstrut@ width\z@}\noexpand\fi%
   \ifStrich@\hbox{\vrule height\htstrut@ depth\dpstrut@ width\z@}%
   \fi$\dsize\m@th\the\hashtoks@$\hfil}%
 \ifStrich@\else\Str@chfehlt#1\\\fi%
 \setbox\z@\hbox{\xdef\Preamble@{\preamble@}}\ifnum`{=0 \fi\iffalse}\fi
 \ialign\bgroup\span\Preamble@\crcr}
\newif\ifhline@\hline@false
\newif\ifhline@@\hline@@false
\def\hlinefor#1{\multispan@{\strip@#1 }\leaders\hrule height\D@cke\hfill%
    \global\hline@true\ignorespaces}
\def\Item "#1"{\par\noindent\hangindent2\parindent%
  \hangafter1\setbox0\hbox{\rm#1\enspace}\ifdim\wd0>2\parindent%
  \box0\else\hbox to 2\parindent{\rm#1\hfil}\fi\ignorespaces}
\def\ITEM #1"#2"{\par\noindent\hangafter1\hangindent#1%
  \setbox0\hbox{\rm#2\enspace}\ifdim\wd0>#1%
  \box0\else\hbox to 0pt{\rm#2\hss}\hskip#1\fi\ignorespaces}
\def\item"#1"{\par\noindent\hang%
  \setbox0=\hbox{\rm#1\enspace}\ifdim\wd0>\the\parindent%
  \box0\else\hbox to \parindent{\rm#1\hfil}\enspace\fi\ignorespaces}
\let\plainitem@\item
\catcode`\@=13

\catcode`\@=11
\font\tenln    = line10
\font\tenlnw   = linew10

\newskip\Einheit \Einheit=0.5cm
\newcount\xcoord \newcount\ycoord
\newdimen\xdim \newdimen\ydim \newdimen\PfadD@cke \newdimen\Pfadd@cke

%%%%%%%%%%%%%%%%%%%%%%%%%%%%%%%%%%%%%%%%%%%%%%%%%
%LaTeX counters, dimensions, variables for lines%
%%%%%%%%%%%%%%%%%%%%%%%%%%%%%%%%%%%%%%%%%%%%%%%%%
\newcount\@tempcnta
\newcount\@tempcntb

\newdimen\@tempdima
\newdimen\@tempdimb

\newdimen\@wholewidth
\newdimen\@halfwidth

\newcount\@xarg
\newcount\@yarg
\newcount\@yyarg
\newbox\@linechar
\newbox\@tempboxa
\newdimen\@linelen
\newdimen\@clnwd
\newdimen\@clnht

\newif\if@negarg

\def\@whilenoop#1{}
\def\@whiledim#1\do #2{\ifdim #1\relax#2\@iwhiledim{#1\relax#2}\fi}
\def\@iwhiledim#1{\ifdim #1\let\@nextwhile=\@iwhiledim
        \else\let\@nextwhile=\@whilenoop\fi\@nextwhile{#1}}

\def\@whileswnoop#1\fi{}
\def\@whilesw#1\fi#2{#1#2\@iwhilesw{#1#2}\fi\fi}
\def\@iwhilesw#1\fi{#1\let\@nextwhile=\@iwhilesw
         \else\let\@nextwhile=\@whileswnoop\fi\@nextwhile{#1}\fi}

\def\thinlines{\let\@linefnt\tenln \let\@circlefnt\tencirc
  \@wholewidth\fontdimen8\tenln \@halfwidth .5\@wholewidth}
\def\thicklines{\let\@linefnt\tenlnw \let\@circlefnt\tencircw
  \@wholewidth\fontdimen8\tenlnw \@halfwidth .5\@wholewidth}
\thinlines
%%%%%%%%%%%%%%%%%%%%%%%%%%%%%%%%%%%%%%%%%%%%%%%%%%%%%%%%%%%

\PfadD@cke1pt \Pfadd@cke0.5pt
\def\PfadDicke#1{\PfadD@cke#1 \divide\PfadD@cke by2 \Pfadd@cke\PfadD@cke \multiply\PfadD@cke by2}
\long\def\LOOP#1\REPEAT{\def\BODY{#1}\ITERATE}
\def\ITERATE{\BODY \let\next\ITERATE \else\let\next\relax\fi \next}
\let\REPEAT=\fi
\def\Punkt{\hbox{\raise-2pt\hbox to0pt{\hss$\ssize\bullet$\hss}}}
\def\DuennPunkt(#1,#2){\unskip
  \raise#2 \Einheit\hbox to0pt{\hskip#1 \Einheit
          \raise-2.5pt\hbox to0pt{\hss$\bullet$\hss}\hss}}
\def\NormalPunkt(#1,#2){\unskip
  \raise#2 \Einheit\hbox to0pt{\hskip#1 \Einheit
          \raise-3pt\hbox to0pt{\hss\twelvepoint$\bullet$\hss}\hss}}
\def\DickPunkt(#1,#2){\unskip
  \raise#2 \Einheit\hbox to0pt{\hskip#1 \Einheit
          \raise-4pt\hbox to0pt{\hss\fourteenpoint$\bullet$\hss}\hss}}
\def\Kreis(#1,#2){\unskip
  \raise#2 \Einheit\hbox to0pt{\hskip#1 \Einheit
          \raise-4pt\hbox to0pt{\hss\fourteenpoint$\circ$\hss}\hss}}

%%%%%%%%%%%%%%%%%%%%%
%LaTeX line macros%
%%%%%%%%%%%%%%%%%%%%%
\def\Line@(#1,#2)#3{\@xarg #1\relax \@yarg #2\relax
\@linelen=#3\Einheit
\ifnum\@xarg =0 \@vline
  \else \ifnum\@yarg =0 \@hline \else \@sline\fi
\fi}

\def\@sline{\ifnum\@xarg< 0 \@negargtrue \@xarg -\@xarg \@yyarg -\@yarg
  \else \@negargfalse \@yyarg \@yarg \fi
\ifnum \@yyarg >0 \@tempcnta\@yyarg \else \@tempcnta -\@yyarg \fi
\ifnum\@tempcnta>6 \@badlinearg\@tempcnta0 \fi
\ifnum\@xarg>6 \@badlinearg\@xarg 1 \fi
\setbox\@linechar\hbox{\@linefnt\@getlinechar(\@xarg,\@yyarg)}%
\ifnum \@yarg >0 \let\@upordown\raise \@clnht\z@
   \else\let\@upordown\lower \@clnht \ht\@linechar\fi
\@clnwd=\wd\@linechar
\if@negarg \hskip -\wd\@linechar \def\@tempa{\hskip -2\wd\@linechar}\else
     \let\@tempa\relax \fi
\@whiledim \@clnwd <\@linelen \do
  {\@upordown\@clnht\copy\@linechar
   \@tempa
   \advance\@clnht \ht\@linechar
   \advance\@clnwd \wd\@linechar}%
\advance\@clnht -\ht\@linechar
\advance\@clnwd -\wd\@linechar
\@tempdima\@linelen\advance\@tempdima -\@clnwd
\@tempdimb\@tempdima\advance\@tempdimb -\wd\@linechar
\if@negarg \hskip -\@tempdimb \else \hskip \@tempdimb \fi
\multiply\@tempdima \@m
\@tempcnta \@tempdima \@tempdima \wd\@linechar \divide\@tempcnta \@tempdima
\@tempdima \ht\@linechar \multiply\@tempdima \@tempcnta
\divide\@tempdima \@m
\advance\@clnht \@tempdima
\ifdim \@linelen <\wd\@linechar
   \hskip \wd\@linechar
  \else\@upordown\@clnht\copy\@linechar\fi}

\def\@hline{\ifnum \@xarg <0 \hskip -\@linelen \fi
\vrule height\Pfadd@cke width \@linelen depth\Pfadd@cke
\ifnum \@xarg <0 \hskip -\@linelen \fi}

\def\@getlinechar(#1,#2){\@tempcnta#1\relax\multiply\@tempcnta 8
\advance\@tempcnta -9 \ifnum #2>0 \advance\@tempcnta #2\relax\else
\advance\@tempcnta -#2\relax\advance\@tempcnta 64 \fi
\char\@tempcnta}

\def\Vektor(#1,#2)#3(#4,#5){\unskip\leavevmode
  \xcoord#4\relax \ycoord#5\relax
      \raise\ycoord \Einheit\hbox to0pt{\hskip\xcoord \Einheit
         \Vector@(#1,#2){#3}\hss}}

\def\Vector@(#1,#2)#3{\@xarg #1\relax \@yarg #2\relax
\@tempcnta \ifnum\@xarg<0 -\@xarg\else\@xarg\fi
\ifnum\@tempcnta<5\relax
\@linelen=#3\Einheit
\ifnum\@xarg =0 \@vvector
  \else \ifnum\@yarg =0 \@hvector \else \@svector\fi
\fi
\else\@badlinearg\fi}

\def\@hvector{\@hline\hbox to 0pt{\@linefnt
\ifnum \@xarg <0 \@getlarrow(1,0)\hss\else
    \hss\@getrarrow(1,0)\fi}}

\def\@vvector{\ifnum \@yarg <0 \@downvector \else \@upvector \fi}

\def\@svector{\@sline
\@tempcnta\@yarg \ifnum\@tempcnta <0 \@tempcnta=-\@tempcnta\fi
\ifnum\@tempcnta <5
  \hskip -\wd\@linechar
  \@upordown\@clnht \hbox{\@linefnt  \if@negarg
  \@getlarrow(\@xarg,\@yyarg) \else \@getrarrow(\@xarg,\@yyarg) \fi}%
\else\@badlinearg\fi}

\def\@upline{\hbox to \z@{\hskip -.5\Pfadd@cke \vrule width \Pfadd@cke
   height \@linelen depth \z@\hss}}

\def\@downline{\hbox to \z@{\hskip -.5\Pfadd@cke \vrule width \Pfadd@cke
   height \z@ depth \@linelen \hss}}

\def\@upvector{\@upline\setbox\@tempboxa\hbox{\@linefnt\char'66}\raise
     \@linelen \hbox to\z@{\lower \ht\@tempboxa\box\@tempboxa\hss}}

\def\@downvector{\@downline\lower \@linelen
      \hbox to \z@{\@linefnt\char'77\hss}}

\def\@getlarrow(#1,#2){\ifnum #2 =\z@ \@tempcnta='33\else
\@tempcnta=#1\relax\multiply\@tempcnta \sixt@@n \advance\@tempcnta
-9 \@tempcntb=#2\relax\multiply\@tempcntb \tw@
\ifnum \@tempcntb >0 \advance\@tempcnta \@tempcntb\relax
\else\advance\@tempcnta -\@tempcntb\advance\@tempcnta 64
\fi\fi\char\@tempcnta}

\def\@getrarrow(#1,#2){\@tempcntb=#2\relax
\ifnum\@tempcntb < 0 \@tempcntb=-\@tempcntb\relax\fi
\ifcase \@tempcntb\relax \@tempcnta='55 \or
\ifnum #1<3 \@tempcnta=#1\relax\multiply\@tempcnta
24 \advance\@tempcnta -6 \else \ifnum #1=3 \@tempcnta=49
\else\@tempcnta=58 \fi\fi\or
\ifnum #1<3 \@tempcnta=#1\relax\multiply\@tempcnta
24 \advance\@tempcnta -3 \else \@tempcnta=51\fi\or
\@tempcnta=#1\relax\multiply\@tempcnta
\sixt@@n \advance\@tempcnta -\tw@ \else
\@tempcnta=#1\relax\multiply\@tempcnta
\sixt@@n \advance\@tempcnta 7 \fi\ifnum #2<0 \advance\@tempcnta 64 \fi
\char\@tempcnta}
%%%%%%%%%%%%%%%%%%%%%%%%%%%%%%%%%%%%%%%%%%%%%%%%%%%%%%%%%%%%%

\def\Diagonale(#1,#2)#3{\unskip\leavevmode
  \xcoord#1\relax \ycoord#2\relax
      \raise\ycoord \Einheit\hbox to0pt{\hskip\xcoord \Einheit
         \Line@(1,1){#3}\hss}}
\def\AntiDiagonale(#1,#2)#3{\unskip\leavevmode
  \xcoord#1\relax \ycoord#2\relax %\advance\xcoord by -0.05\relax
      \raise\ycoord \Einheit\hbox to0pt{\hskip\xcoord \Einheit
         \Line@(1,-1){#3}\hss}}
\def\Pfad(#1,#2),#3\endPfad{\unskip\leavevmode
  \xcoord#1 \ycoord#2 \thicklines\ZeichnePfad#3\endPfad\thinlines}
\def\ZeichnePfad#1{\ifx#1\endPfad\let\next\relax
  \else\let\next\ZeichnePfad
    \ifnum#1=1
      \raise\ycoord \Einheit\hbox to0pt{\hskip\xcoord \Einheit
         \vrule height\Pfadd@cke width1 \Einheit depth\Pfadd@cke\hss}%
      \advance\xcoord by 1
    \else\ifnum#1=2
      \raise\ycoord \Einheit\hbox to0pt{\hskip\xcoord \Einheit
        \hbox{\hskip-\PfadD@cke\vrule height1 \Einheit width\PfadD@cke depth0pt}\hss}%
      \advance\ycoord by 1
    \else\ifnum#1=3
      \raise\ycoord \Einheit\hbox to0pt{\hskip\xcoord \Einheit
         \Line@(1,1){1}\hss}
      \advance\xcoord by 1
      \advance\ycoord by 1
    \else\ifnum#1=4
      \raise\ycoord \Einheit\hbox to0pt{\hskip\xcoord \Einheit
         \Line@(1,-1){1}\hss}
      \advance\xcoord by 1
      \advance\ycoord by -1
    \else\ifnum#1=5
      \advance\xcoord by -1
      \raise\ycoord \Einheit\hbox to0pt{\hskip\xcoord \Einheit
         \vrule height\Pfadd@cke width1 \Einheit depth\Pfadd@cke\hss}%
    \else\ifnum#1=6
      \advance\ycoord by -1
      \raise\ycoord \Einheit\hbox to0pt{\hskip\xcoord \Einheit
        \hbox{\hskip-\PfadD@cke\vrule height1 \Einheit width\PfadD@cke depth0pt}\hss}%
    \else\ifnum#1=7
      \advance\xcoord by -1
      \advance\ycoord by -1
      \raise\ycoord \Einheit\hbox to0pt{\hskip\xcoord \Einheit
         \Line@(1,1){1}\hss}
    \else\ifnum#1=8
      \advance\xcoord by -1
      \advance\ycoord by +1
      \raise\ycoord \Einheit\hbox to0pt{\hskip\xcoord \Einheit
         \Line@(1,-1){1}\hss}
    \fi\fi\fi\fi
    \fi\fi\fi\fi
  \fi\next}
\def\hSSchritt{\leavevmode\raise-.4pt\hbox to0pt{\hss.\hss}\hskip.2\Einheit
  \raise-.4pt\hbox to0pt{\hss.\hss}\hskip.2\Einheit
  \raise-.4pt\hbox to0pt{\hss.\hss}\hskip.2\Einheit
  \raise-.4pt\hbox to0pt{\hss.\hss}\hskip.2\Einheit
  \raise-.4pt\hbox to0pt{\hss.\hss}\hskip.2\Einheit}
\def\vSSchritt{\vbox{\baselineskip.2\Einheit\lineskiplimit0pt
\hbox{.}\hbox{.}\hbox{.}\hbox{.}\hbox{.}}}
\def\DSSchritt{\leavevmode\raise-.4pt\hbox to0pt{%
  \hbox to0pt{\hss.\hss}\hskip.2\Einheit
  \raise.2\Einheit\hbox to0pt{\hss.\hss}\hskip.2\Einheit
  \raise.4\Einheit\hbox to0pt{\hss.\hss}\hskip.2\Einheit
  \raise.6\Einheit\hbox to0pt{\hss.\hss}\hskip.2\Einheit
  \raise.8\Einheit\hbox to0pt{\hss.\hss}\hss}}
\def\dSSchritt{\leavevmode\raise-.4pt\hbox to0pt{%
  \hbox to0pt{\hss.\hss}\hskip.2\Einheit
  \raise-.2\Einheit\hbox to0pt{\hss.\hss}\hskip.2\Einheit
  \raise-.4\Einheit\hbox to0pt{\hss.\hss}\hskip.2\Einheit
  \raise-.6\Einheit\hbox to0pt{\hss.\hss}\hskip.2\Einheit
  \raise-.8\Einheit\hbox to0pt{\hss.\hss}\hss}}
\def\SPfad(#1,#2),#3\endSPfad{\unskip\leavevmode
  \xcoord#1 \ycoord#2 \ZeichneSPfad#3\endSPfad}
\def\ZeichneSPfad#1{\ifx#1\endSPfad\let\next\relax
  \else\let\next\ZeichneSPfad
    \ifnum#1=1
      \raise\ycoord \Einheit\hbox to0pt{\hskip\xcoord \Einheit
         \hSSchritt\hss}%
      \advance\xcoord by 1
    \else\ifnum#1=2
      \raise\ycoord \Einheit\hbox to0pt{\hskip\xcoord \Einheit
        \hbox{\hskip-2pt \vSSchritt}\hss}%
      \advance\ycoord by 1
    \else\ifnum#1=3
      \raise\ycoord \Einheit\hbox to0pt{\hskip\xcoord \Einheit
         \DSSchritt\hss}
      \advance\xcoord by 1
      \advance\ycoord by 1
    \else\ifnum#1=4
      \raise\ycoord \Einheit\hbox to0pt{\hskip\xcoord \Einheit
         \dSSchritt\hss}
      \advance\xcoord by 1
      \advance\ycoord by -1
    \else\ifnum#1=5
      \advance\xcoord by -1
      \raise\ycoord \Einheit\hbox to0pt{\hskip\xcoord \Einheit
         \hSSchritt\hss}%
    \else\ifnum#1=6
      \advance\ycoord by -1
      \raise\ycoord \Einheit\hbox to0pt{\hskip\xcoord \Einheit
        \hbox{\hskip-2pt \vSSchritt}\hss}%
    \else\ifnum#1=7
      \advance\xcoord by -1
      \advance\ycoord by -1
      \raise\ycoord \Einheit\hbox to0pt{\hskip\xcoord \Einheit
         \DSSchritt\hss}
    \else\ifnum#1=8
      \advance\xcoord by -1
      \advance\ycoord by 1
      \raise\ycoord \Einheit\hbox to0pt{\hskip\xcoord \Einheit
         \dSSchritt\hss}
    \fi\fi\fi\fi
    \fi\fi\fi\fi
  \fi\next}
\def\Koordinatenachsen(#1,#2){\unskip
 \hbox to0pt{\hskip-.5pt\vrule height#2 \Einheit width.5pt depth1 \Einheit}%
 \hbox to0pt{\hskip-1 \Einheit \xcoord#1 \advance\xcoord by1
    \vrule height0.25pt width\xcoord \Einheit depth0.25pt\hss}}
\def\Koordinatenachsen(#1,#2)(#3,#4){\unskip
 \hbox to0pt{\hskip-.5pt \ycoord-#4 \advance\ycoord by1
    \vrule height#2 \Einheit width.5pt depth\ycoord \Einheit}%
 \hbox to0pt{\hskip-1 \Einheit \hskip#3\Einheit 
    \xcoord#1 \advance\xcoord by1 \advance\xcoord by-#3 
    \vrule height0.25pt width\xcoord \Einheit depth0.25pt\hss}}
\def\Gitter(#1,#2){\unskip \xcoord0 \ycoord0 \leavevmode
  \LOOP\ifnum\ycoord<#2
    \loop\ifnum\xcoord<#1
      \raise\ycoord \Einheit\hbox to0pt{\hskip\xcoord \Einheit\Punkt\hss}%
      \advance\xcoord by1
    \repeat
    \xcoord0
    \advance\ycoord by1
  \REPEAT}
\def\Gitter(#1,#2)(#3,#4){\unskip \xcoord#3 \ycoord#4 \leavevmode
  \LOOP\ifnum\ycoord<#2
    \loop\ifnum\xcoord<#1
      \raise\ycoord \Einheit\hbox to0pt{\hskip\xcoord \Einheit\Punkt\hss}%
      \advance\xcoord by1
    \repeat
    \xcoord#3
    \advance\ycoord by1
  \REPEAT}
\def\Label#1#2(#3,#4){\unskip \xdim#3 \Einheit \ydim#4 \Einheit
  \def\lo{\advance\xdim by-.5 \Einheit \advance\ydim by.5 \Einheit}%
  \def\llo{\advance\xdim by-.25cm \advance\ydim by.5 \Einheit}%
  \def\loo{\advance\xdim by-.5 \Einheit \advance\ydim by.25cm}%
  \def\o{\advance\ydim by.25cm}%
  \def\ro{\advance\xdim by.5 \Einheit \advance\ydim by.5 \Einheit}%
  \def\rro{\advance\xdim by.25cm \advance\ydim by.5 \Einheit}%
  \def\roo{\advance\xdim by.5 \Einheit \advance\ydim by.25cm}%
  \def\l{\advance\xdim by-.30cm}%
  \def\r{\advance\xdim by.30cm}%
  \def\lu{\advance\xdim by-.5 \Einheit \advance\ydim by-.6 \Einheit}%
  \def\llu{\advance\xdim by-.25cm \advance\ydim by-.6 \Einheit}%
  \def\luu{\advance\xdim by-.5 \Einheit \advance\ydim by-.30cm}%
  \def\u{\advance\ydim by-.30cm}%
  \def\ru{\advance\xdim by.5 \Einheit \advance\ydim by-.6 \Einheit}%
  \def\rru{\advance\xdim by.25cm \advance\ydim by-.6 \Einheit}%
  \def\ruu{\advance\xdim by.5 \Einheit \advance\ydim by-.30cm}%
  #1\raise\ydim\hbox to0pt{\hskip\xdim
     \vbox to0pt{\vss\hbox to0pt{\hss$#2$\hss}\vss}\hss}%
}
\catcode`\@=13

\hsize13cm
\vsize19cm
\newdimen\fullhsize
\newdimen\fullvsize
\newdimen\halfsize
\fullhsize13cm
\fullvsize19cm
\halfsize=0.5\fullhsize
\advance\halfsize by-0.5em

\magnification1200

\TagsOnRight

\def\AignAB{1}
\def\CaYiAA{2}
\def\CiglAS{3}
\def\CiglAT{4}
\def\CiglAU{5}
\def\CvRIAA{6}
\def\DoFSAA{7}
\def\ElouAA{8}
\def\EuWYAA{9}
\def\GaRaAF{10}
\def\GeViAA{11}
\def\GeViAB{12}
\def\HaReAA{13}
\def\IsStAB{14}
\def\IsStAC{15}
\def\KratBN{16}
\def\KratBZ{17}
\def\KratCI{18}
\def\KratCO{19}
\def\LindAA{20}
\def\MuWYAA{21}
\def\OEIS{22}
\def\StanBI{23}
\def\VienAE{24}

%equation numbers
\def\AA{1.1}
\def\AB{1.2}
\def\AC{1.3}
\def\AL{1.4}
\def\BA{2.1}
\def\BB{2.2}
\def\BBa{2.3}
\def\BE{2.4}
\def\BC{2.5}
\def\BD{2.6}
\def\AD{2.7}
\def\AE{2.8}
\def\AI{2.9}
\def\AJ{2.10}
\def\AF{2.11}
\def\AG{2.12}
\def\CA{3.1}
\def\CB{3.2}
\def\EA{3.3}
\def\EB{3.4}
\def\EC{3.5}
\def\ED{3.6}
\def\EF{4.1}
\def\EFa{4.2}
\def\EG{4.3}
\def\EH{4.4}
\def\EHa{4.5}
\def\YA{4.6}
\def\YB{4.7}
\def\YC{4.8}
\def\YD{4.9}
\def\YE{4.10}
\def\YF{4.11}
\def\YG{4.12}
\def\YGa{4.13}
\def\YGb{4.14}
\def\AX{4.15}
\def\AY{4.16}
\def\YH{4.17}
\def\YI{4.18}
\def\YIa{4.19}
\def\YIb{4.20}
\def\YIc{4.21}
\def\YJ{4.22}
\def\YK{4.23}
\def\YL{4.24}
\def\YO{4.25}
\def\YP{4.26}
\def\YQ{4.27}
\def\YR{4.28}
\def\YS{4.29}
\def\AU{4.30}
\def\AV{4.31}
\def\YM{4.32}
\def\YN{4.33}
\def\AK{4.34}
\def\AM{4.35}
\def\AMa{4.36}
\def\AMb{4.37}
\def\AO{4.38}
\def\AR{4.39}
\def\AS{4.40}

\def\AN{5.1}
\def\CC{5.2}
\def\CD{5.3}
\def\CE{5.4}
\def\CF{5.5}
\def\CG{5.6}
\def\CH{5.7}

\def\DA{6.1}
\def\DD{6.2}
\def\AW{6.3}
\def\DB{6.4}
\def\DC{6.5}

\def\ZA{7.1}
\def\ZB{7.2}
\def\ZC{7.3}
\def\ZD{7.4}
\def\ZE{7.5}
\def\ZF{7.6}
\def\ZG{7.7}
\def\ZH{7.8}
\def\ZI{7.9}

%theorem numbers
%\def\TA{1}
%\def\TB{2}
\def\UA{1}
\def\TC{2}
\def\TG{3}
\def\TE{4}
\def\TD{5}
\def\TF{6}
\def\TI{7}
\def\TH{8}

%

%figure numbers
\def\FA{1}
\def\FB{2}
\def\FC{3}
\def\FD{4}
\def\FE{5}
\def\FF{6}
\def\FG{7}

\def\GF{\operatorname{GF}}

\def\fl#1{\lfloor#1\rfloor}

\topmatter 
\title Hankel determinants of linear combinations
of moments of orthogonal polynomials
\endtitle 
\author J. Cigler and C.~Krattenthaler$^{\dagger}$
\endauthor 
\affil 
Fakult\"at f\"ur Mathematik, Universit\"at Wien,\\
Oskar-Morgenstern-Platz~1, A-1090 Vienna, Austria.\\
WWW: {\tt http://homepage.univie.ac.at/johann.cigler}\\
WWW: \tt http://www.mat.univie.ac.at/\~{}kratt
\endaffil
\address Fakult\"at f\"ur Mathematik, Universit\"at Wien,
Oskar-Morgenstern-Platz~1, A-1090 Vienna, Austria.\newline
WWW: \tt http://homepage.univie.ac.at/johann.cigler,
http://www.mat.univie.ac.at/\~{}kratt
\endaddress
\thanks $^\dagger$Research partially supported by the Austrian
Science Foundation FWF (grant S50-N15)
in the framework of the Special Research Program
``Algorithmic and Enumerative Combinatorics"%
\endthanks

\subjclass Primary 05A19;
 Secondary 05A10 05A15 11C20 15A15 33C45 42C05
\endsubjclass
\keywords Hankel determinants, moments of orthogonal polynomials,
Catalan numbers, 
Motzkin numbers, Fibonacci numbers, Lucas numbers,
Schr\"oder numbers, Riordan numbers, Fine numbers,
central binomial coefficients, central trinomial numbers,
Chebyshev polynomials, Rogers--Szeg\H o polynomials,
Dyck paths, Motzkin paths, 
non-intersecting lattice paths
\endkeywords
\abstract 
We prove evaluations of 
Hankel determinants of linear combinations of moments of
orthogonal polynomials (or, equivalently, of generating functions for
Motzkin paths), thus generalising known results for Catalan numbers.
\endabstract
\endtopmatter
\document

\subhead 1. Introduction\endsubhead
Let $C_n=\frac {1} {n+1}\binom {2n}n$ denote the $n$-th {\it Catalan number}, 
and let $(F_n)_{n\ge0}$ be the sequence of {\it Fibonacci
numbers} defined by $F_n=F_{n-1}+F_{n-2}$ with initial values 
$F_0=0$ and $F_1=1$. 

In \cite{\CvRIAA}, Cvetkovi\'c, Rajkovi\'c and Ivkovi\'c proved the
Hankel determinant evaluations
$$
\det\left(C_{i+j}+C_{i+j+1}\right)_{i,j=0}^{n-1}
=F_{2n+1}
\tag\AA
$$
and
$$
\det\left(C_{i+j+1}+C_{i+j+2}\right)_{i,j=0}^{n-1}
=F_{2n+2}.
\tag\AB
$$

In \cite{\DoFSAA}, Dougherty, French, Saderholm and Qian proved
a result of similar flavour, namely
$$
\det\left(C_{i+j}+2C_{i+j+1}+C_{i+j+2}\right)_{i,j=0}^{n-1}
=\sum_{j=0}^n F_{2j+1}^2.
\tag\AC
$$
One finds variations and generalisations of these identities at
numerous places in the literature, for example in
\cite{\CaYiAA, \CiglAS, \CiglAU, \CvRIAA, \DoFSAA, \ElouAA, \EuWYAA,
\KratCI, \MuWYAA}.

Very recently, the first author became aware of
$$
\frac{\det\left(\binom {2i+2j+2}{i+j+1}+2\binom {2i+2j+4}{i+j+2}
+\binom {2i+2j+6}{i+j+3}\right)_{i,j=0}^{n-1}}
{2^{n}}
=\sum_{j=0}^n L_{2j+1}^2,
\tag\AL
$$
where $L_n$ denotes the $n$-th {\it Lucas number}, defined by the
Fibonacci recurrence $L_n=L_{n-1}+L_{n-2}$ with initial values
$L_0=2$ and $L_1=1$.
This formula was posted without proof by a mathematical amateur, Tony Foster,
in a Facebook group about Pascal's triangle.

We decided to search for the general background of this kind of Hankel
determinant evaluations. Our main result in Theorem~\UA\ 
in the next section
provides a formula for the evaluation of Hankel determinants of
a linear combination of three generating functions for Motzkin paths.
In these generating functions, the weight of a path depends on its steps.
In that section, in Corollary~\TC, we then record separately
the corresponding result for Dyck paths. 
Indeed, by appropriate specialisations,
all the above identities as well as many more, in the literature or
not, can be obtained, see Sections~3 and~4. One result that is of
particular note is Corollary~\TG, which gives a closed form
evaluation, in terms of Chebyshev polynomials, 
of Hankel determinants as in Theorem~\UA\ for the
special case where weights of steps do not depend on their height,
except if they are at height zero.

We provide two different proofs for Theorem~\UA\ and Corollary~\TC,
both of which are of interest as we believe. The proofs in Section~5
are based on the decomposition of the matrices of which we want to
compute the determinant into the product of ``easier to handle''
matrices. In particular, the proof of Corollary~\TC\ in that section
shows that it indeed follows from our main theorem, Theorem~\UA,
although this is not completely obvious from the outset.
On the other hand, in Section~6, we present combinatorial
proofs that make use of the combinatorics of non-intersecting lattice
paths. Once the setup is explained, the proofs themselves then
consist of just one picture in each case. 
We point out that the basic
idea has appeared before in \cite{\KratBZ, paragraph after Theorem~29
and paragraph above Theorem~30} and
\cite{\CaYiAA} but has not been fully exploited there.

Finally, in Section~7 we consider the Hankel determinants of linear
combinations of {\it four} path generating functions. Indeed, 
the combinatorial model from Section~6 can be readily extended to
yield a corresponding formula, see Theorem~\TD. 
We also derive the extension of Corollary~\TG\
for the specialisation where the step weights do not depend on the
height of the steps (except if they are at height zero), see Corollary~\TF. 
By a perusal of Corollaries~\TG\ and~\TF, a pattern emerges.
On this basis, we make a conjecture on the evaluation
of the Hankel determinant of an {\it arbitrarily long} linear combination of
path generating functions in the case of the aforementioned
specialisation, see Conjecture~\TH.

\subhead 2. The main theorem and a corollary\endsubhead
Let $(s_n)_{n\ge0}$ and $(t_n)_{n\ge0}$
be sequences of real numbers with $t_n\ne0$ for all $n$, and define
polynomials (in $\al$ and the $s_i$'s and $t_i$'s) $f_n(\al)$ by the recursion
$$
f_n(\al)=(\al+s_{n-1})f_{n-1}(\al)-t_{n-2}f_{n-2}(\al),
\tag\BA
$$
with $f_0(\al) = 1$ and $f_{-1}(\al)=0$.

Consider now {\it Motzkin paths}, i.e., lattice paths on
the two-dimensional integer lattice $\Z^2$ with up-steps
$(x,y)\to(x +1,y +1)$, horizontal steps $(x,y)\to(x +1,y)$, and down-steps
$(x,y)\to(x +1,y -1)$, which start at
$(0, 0)$ and never run below the $x$-axis.\footnote{We should point out
that we use a slightly extended version of the notion of ``Motzkin path'' here.
Usually, a Motzkin path is required to return to the $x$-axis at the
end. We do not make this requirement.} The weight of such
a path is the
product of the weights of its steps, where the weight of an up-step is
1, the weight of a
horizontal step on height $h$ is $s_h$, and the weight of a down-step 
which ends on height $h$ is $t_h$.

Let $m(n,k)$ denote the sum of the weights of all such paths with end point
$(n,k)$. Then we have
$$
m(n,k)=m(n-1,k-1)+s_km(n-1,k)+t_km(n-1,k+1),
\tag\BB
$$
with $m(0,k) = [k = 0]$ and $m(n,k) = 0$ for $k < 0$,
where the Iverson bracket $[\Cal A]$ is defined by
$[\Cal A]=1$ if $\Cal A$ is true and $[\Cal A]=0$ otherwise.
In order to explain the title of the paper:
it is well known (see e.g\. \cite{\KratBN, Theorems~11--13}) 
that the polynomials $p_n(x)$ defined recursively by 
$$p_n(x)=(x-s_{n-1})p_{n-1}(x)-t_{n-2}p_{n-2}(x),
\tag\BBa
$$
with initial values
$p_{-1} (x)= 0$ and $p_0 (x)=1$ are orthogonal with respect to the linear
functional $L$ defined by $L(p_n(x))=[n=0]$,
and that their moments are $L(x^n)=m(n,0)$.
(In other words: the polynomials $f_n(\al)$ are related to the
polynomials $p_n(\al)$ by the replacement of $s_i$ by $-s_i$ for all~$i$.)
So let us write $m_n$ for $m(n,0)$. The reader should be aware that,
if we set $s_i=t_i=1$ for all~$i$, then $m_n$ reduces to the classical
{\it Motzkin number} $M_n$ given by 
$M_n=\sum _{k=0} ^{\fl{n/2}}\binom {n} {2k} \frac {1} {k+1}\binom
{2k}k$, 
and $f_n(\al)$ becomes $U_n\big((\al+1)/2\big)$, where $U_n(x)$
is the $n$-th {\it Chebyshev polynomial of the second kind\/}
$$
U_n(x)=\sum_{k\ge0}(-1)^k\binom {n-k}{k}(2x)^{n-2k}.
\tag\BE
$$

With these notations we are ready to state our main result
about Hankel determinants of linear combinations of the path
generating functions $m_n=m(n,0)$.

\proclaim{Theorem \UA}
Let $\al$ and $\be$ be indeterminates. Then,
for all positive integers $n$, we have
$$
\frac {\det\left(\al\be m_{i+j}+(\al+\be)m_{i+j+1}+m_{i+j+2}\right)_{i,j=0}^{n-1}} 
{\det\left(m_{i+j}\right)_{i,j=0}^{n-1}}
=\sum_{j=0}^nf_{j}(\al)f_j(\be)
\prod _{\ell=j} ^{n-1}t_\ell.
\tag\BC
$$
\endproclaim

\remark{Remark}
We should point out that it is well known (see e.g\. \cite{\VienAE,
  Ch.~IV, Cor.~6}) that
$$
{\det\left(m_{i+j}\right)_{i,j=0}^{n-1}}
=\prod _{i=0} ^{n-1}t_i^{n-i-1}.
\tag\BD
$$
\endremark

\medskip
It is worthwhile to state the analogous result for {\it Dyck paths}
separately.

Let $(T_n)_{n\ge0}$ be a sequence of real numbers with $T_n\ne0$ for all~$n$.
Define polynomials (in $\al$ and the $T_i$'s) $g_n(\al)$ by
$$
\align
g_{2n}(\al)&=\al g_{2n-1}(\al)+T_{2n-2}g_{2n-2}(\al),\\
g_{2n+1}(\al)&=g_{2n}(\al)+T_{2n-1}g_{2n-1}(\al),
\tag\AD
\endalign
$$
with $g_{-1}(\al)=0$ and $g_0(\al)=1$.

Consider now Dyck paths, i.e., lattice paths on
the two-dimensional integer lattice $\Z^2$ with up-steps
$(x,y)\to(x +1,y +1)$ and down-steps
$(x,y)\to(x +1,y -1)$, which  start at
$(0, 0)$ and never run below the $x$-axis.\footnote{Again, we are
using a slightly extended version of the notion of ``Dyck path'' here,
by {\it not\/} requiring a Dyck path to return to the $x$-axis at the end.} 
As before, the weight of such a path is the
product of the weights of its steps. Here, the weight of an up-step is
1, and the weight of a down-step 
which ends on height $h$ is $T_h$.

Let $c(n,k)$ denote the sum of the weights of all such paths with end point
$(n,k)$. Then we have
$$
c(n,k)=c(n-1,k-1)+T_kc(n-1,k+1),
\tag\AE
$$
with $c(0,k) =[k = 0]$ and $c(n,-1) = 0$.
Finally, set  
$c_n= c(2n, 0)$.
Here, the reader should be aware that,
if we set $T_i=1$ for all~$i$, then $c_n$ reduces to the Catalan
number $C_n$. Furthermore, it can be shown that 
$g_{2n}(\al)=U_n((\al+2)/2)-U_{n-1}((\al+2)/2)$ and
$g_{2n+1}(\al)=U_n((\al+2)/2)$,
where $U_n(x)$ denotes as before the $n$-th Chebyshev polynomial~(\BE).

\proclaim{Corollary \TC}
Let $\al$ and $\be$ be indeterminates. Then,
for all positive integers $n$, we have
$$\align 
\frac {\det\left(\al\be c_{i+j}+(\al+\be)c_{i+j+1}+c_{i+j+2}\right)_{i,j=0}^{n-1}} 
{\det\left(c_{i+j}\right)_{i,j=0}^{n-1}}
&=\sum_{j=0}^nT_{2j}T_{2j+1}\cdots T_{2n-1}g_{2j}(\al)g_{2j}(\be),
\tag\AI\\
\frac {\det\left(\al\be c_{i+j+1}+(\al+\be)c_{i+j+2}+c_{i+j+3}\right)_{i,j=0}^{n-1}} 
{\det\left(c_{i+j+1}\right)_{i,j=0}^{n-1}}
&=\sum_{j=0}^nT_{2j+1}T_{2j+2}\cdots T_{2n}g_{2j+1}(\al)g_{2j+1}(\be).
\tag\AJ
\endalign$$
\endproclaim

\remark{Remark}
Again, we must point out that it is well known (see e.g\. \cite{\KratBZ,
  Eqs.~(5.44) and (5.45)}) that
$$
{\det\left(c_{i+j}\right)_{i,j=0}^{n-1}}
=\prod _{i=0} ^{n-1}(T_{2i}T_{2i+1})^{n-i-1}.
\tag\AF
$$
and
$$
{\det\left(c_{i+j+1}\right)_{i,j=0}^{n-1}}
=\prod _{i=0} ^{n-1}T_0(T_{2i+1}T_{2i+2})^{n-i-1}.
\tag\AG
$$
\endremark

\subhead 3. Immediate consequences\endsubhead
In this section, we state some direct implications
of Theorem~\UA\ and of Corollary~\TC\ explicitly which result
by specialising $\be$  to~$0$ respectively to~$\infty$.

First, if we let $\be=0$ in (\BC), then we obtain
$$
\frac {\det\left(\al m_{i+j+1}+m_{i+j+2}\right)_{i,j=0}^{n-1}} 
{\det\left(m_{i+j}\right)_{i,j=0}^{n-1}}
=\sum_{j=0}^nf_{j}(\al)f_j(0)
\prod _{\ell=j} ^{n-1}t_\ell.
\tag\CA
$$

Next we read coefficients of $\be^n$ on both sides of (\BC).
Since each element of the matrix in the numerator 
on the left-hand side of which the determinant is taken is linear
in $\be$ and the matrix itself is an $n\times n$ matrix,
reading the coefficient of $\be^n$ on the left-hand side amounts
to reading the coefficient of $\be^1$ in each element of the matrix.
On the other hand, due to the recurrence (\BA), the coefficient
of $\be^n$ in $f_j(\be)$ is zero as long as $j<n$, whereas it is
$1$ for $f_n(\be)$. Hence, we obtain
$$
\frac {\det\left(\al m_{i+j}+m_{i+j+1}\right)_{i,j=0}^{n-1}} 
{\det\left(m_{i+j}\right)_{i,j=0}^{n-1}}
=f_{n}(\al).
\tag\CB
$$

\remark{Remark}
The identities (\CA) and (\CB) have been given earlier
by Mu, Wang and Yeh in \cite{\MuWYAA, Theorem~1.3}
in a different but equivalent form.
\endremark

In the same way, from Corollary~\TC\ we obtain
$$\align 
\frac {\det\left(\al c_{i+j+1}+c_{i+j+2}\right)_{i,j=0}^{n-1}} 
{\det\left(c_{i+j}\right)_{i,j=0}^{n-1}}
&=\sum_{j=0}^ng_{2j}(\al)g_{2j}(0)
\prod _{\ell=2j} ^{2n-1}T_\ell,
\tag\EA\\
\frac {\det\left(\al c_{i+j+2}+c_{i+j+3}\right)_{i,j=0}^{n-1}} 
{\det\left(c_{i+j+1}\right)_{i,j=0}^{n-1}}
&=\sum_{j=0}^ng_{2j+1}(\al)g_{2j+1}(0)
\prod _{\ell=2j+1} ^{2n}T_\ell,
\tag\EB\\
\frac {\det\left(\al c_{i+j}+c_{i+j+1}\right)_{i,j=0}^{n-1}} 
{\det\left(c_{i+j}\right)_{i,j=0}^{n-1}}
&=g_{2n}(\al),
\tag\EC\\
\frac {\det\left(\al c_{i+j+1}+c_{i+j+2}\right)_{i,j=0}^{n-1}} 
{\det\left(c_{i+j+1}\right)_{i,j=0}^{n-1}}
&=g_{2n+1}(\al).
\tag\ED
\endalign$$

\subhead 4. Specialisations\endsubhead
By (further) specialising $\al$, $\be$, the $s_i$'s and $t_i$'s,
we obtain numerous Hankel determinant evaluations, covering many of
those that one finds in the literature. In particular, they cover all
of the evaluations that we mentioned in the introduction.

In fact, Mu, Wang and Yeh \cite{\MuWYAA} provide an extensive list
of specialisations of the path generating functions $m_n=m(n,0)$ that
give familiar combinatorial sequences. For the convenience of the
reader we reproduce it here (the codes in parentheses refer to the
sequence number in The On-Line Encyclopedia of Integer Sequences
\cite{\OEIS}):

\roster 
\item"(i)"{\it Motzkin numbers} $M_n$: $s_i=t_i=1$ for all $i$
({\tt A001006}).
\item"(ii)"{\it Catalan numbers} $C_n$: $s_0=1$, $s_i=2$ for $i\ge1$, and
$t_i=1$ for all~$i$ ({\tt A000108}).
\item"(iii)"{\it shifted Catalan numbers} $C_{n+1}$: 
$s_i=2$ and $t_i=1$ for all~$i$.
\item"(iv)"{\it central binomial coefficients} $\binom {2n}n$:
$s_i=2$ for all~$i$, $t_0=2$, and $t_i=1$ for $i\ge1$ ({\tt A000984}).
\item"(v)"{\it central trinomial coefficients} $T_n$:
$s_i=1$ for all~$i$, $t_0=2$, and $t_i=1$ for $i\ge1$ ({\tt A002426}).
\item"(vi)"{\it central Delannoy numbers} $D_n$:
$s_i=3$ for all~$i$, $t_0=4$, and $t_i=2$ for $i\ge1$ ({\tt A001850}).
\item"(vii)"{\it large Schr\"oder numbers} $r_n$:
$s_0=2$, $s_i=3$ for $i\ge1$, and
$t_i=2$ for all~$i$ ({\tt A006318}).
\item"(viii)"{\it little Schr\"oder numbers} $S_n$:
$s_0=1$, $s_i=3$ for $i\ge1$, and
$t_i=2$ for all~$i$ ({\tt A001003}).
\item"(ix)"{\it Fine numbers} $\text{Fine}_n$:
$s_0=0$, $s_i=2$ for $i\ge1$, and
$t_i=1$ for all~$i$ ({\tt A000957}).
\item"(x)"{\it Riordan numbers} $R_n$:
$s_0=0$, $s_i=1$ for $i\ge1$, and
$t_i=1$ for all~$i$ ({\tt A005043}).
\item"(xi)"{\it numbers of restricted hexagonal polyominoes} $H_n$:
$s_i=3$ and $t_i=1$ for all~$i$ ({\tt A002212}).
\item"(xii)"{\it Bell numbers} $B_n$:
$s_i=t_i=i+1$ for all $i$ ({\tt A000110}).
\item"(xiii)"{\it factorials} $n!$:
$s_i=2i+1$ and $t_i=(i+1)^2$ for all $i$.
\endroster

We refer the reader to \cite{\AignAB, \StanBI} for definitions and discussion
of most of these numbers, and to \cite{\HaReAA} for the hexagonal
polyominoes (called {\it polyhexes} in \cite{\HaReAA}).

We add a few more specialisations:

\roster 
\item"(xiv)"{\it shifted and scaled central binomial coefficients} 
$\frac {1} {2}\binom {2n+2}{n+1}$:
$s_0=3$, $s_i=2$ for~$i\ge1$, and $t_i=1$ for all~$i$.
\item"(xv)"{\it odd and even central binomial coefficients} 
$\binom {n}{\fl{n/2}}$:
$s_0=1$, $s_i=0$ for~$i\ge1$, and $t_i=1$ for all~$i$.
\item"(xvi)"{\it Catalan numbers interleaved with zeroes}
$[n\text{ even}]C_{n/2}$:
$s_i=0$ and $t_i=1$ for all $i$.
\item"(xvii)"{\it central binomial coefficients interleaved with zeroes}
$[n\text{ even}]\binom n{n/2}$:
$s_i=0$ for all~$i$, $t_0=2$, and $t_i=1$ for $i\ge1$.
\item"(xviii)"{\it Fine numbers interleaved with zeroes}
$[n=0]+[n\text{ even}]\text{Fine}_{n/2}$:
$s_i=0$ for all~$i$, $t_0=1$, and $t_i=-1$ for $i\ge1$.
\endroster

Mu, Wang and Yeh note that, in most of these specialisations, 
$s_i$ and $t_i$ are constant for $i\ge1$ (namely in (i)--(xi)),
and the same is true for (xiv)--(xviii). 
They then discuss that
special case separately, find generating functions for the
values of $m_n$ and $f_n(\al)$ provided that $s_i\equiv s$ and
$t_i\equiv t$ for $i\ge1$, but --- at the time --- missed to notice
that the specialised polynomials $f_n(\al)$ can be expressed in
terms of Chebyshev polynomials, and consequently did not realise
that, for this specialisation, the sum on the right-hand side of
(\CA) can actually be evaluated. The corresponding
(missed) results are the contents of the following corollary.

\proclaim{Corollary \TG}
With the setup in Section~2, let $s_i\equiv s$ and $t_i\equiv t$ for
$i\ge1$. Then $f_0(\al)=1$ and
$$
\align
f_n(\al)&=t_0t^{(n-2)/2}U_n\left(\tfrac {\al+s} {2\sqrt t}\right)
+\big((t-t_0)(\al+s)-t(s-s_0)\big)
t^{(n-3)/2}U_{n-1}\left(\tfrac {\al+s} {2\sqrt t}\right)
\tag\EF
\\
&=t^{n/2}U_n\left(\tfrac {\al+s} {2\sqrt t}\right)
-t^{(n-1)/2}(s-s_0)U_{n-1}\left(\tfrac {\al+s} {2\sqrt t}\right)\\
&\kern4cm
+t^{(n-2)/2}(t-t_0)U_{n-2}\left(\tfrac {\al+s} {2\sqrt t}\right),
\quad \text{for }n\ge1.
\tag\EFa
\endalign
$$
Moreover, in that case we have
$$
\frac {\det\left(\al\be m_{i+j}+(\al+\be)m_{i+j+1}+m_{i+j+2}\right)_{i,j=0}^{n-1}} 
{\det\left(m_{i+j}\right)_{i,j=0}^{n-1}}
=\frac {\text{\rm Num}(\al,\be)} {\al-\be},
\tag\EG
$$
where
$$
\multline
\text{\rm Num}(\al,\be)=
t^{n+1/2}
\big(U_{n+1}(x)U_n(y)-U_n(x)U_{n+1}(y)\big)\\
-t^{n}(s - s_0) 
\big(U_{n+1}(x)U_{n-1}(y)-U_{n-1}(x)U_{n+1}(y)\big)\\
+t^{n-1/2} \big((s -  s_0)^2 - (t - t_0)\big) 
\big(U_{n}(x)U_{n-1}(y)-U_{n-1}(x)U_{n}(y)\big)\\
+t^{n-1/2} (t - t_0)
\big(U_{n+1}(x)U_{n-2}(y)-U_{n-2}(x)U_{n+1}(y)\big)\\
-t^{n-1} (s - s_0) (t-t_0)
\big(U_{n}(x)U_{n-2}(y)-U_{n-2}(x)U_{n}(y)\big)\\
+t^{n-3/2} (t - t_0)^2
\big(U_{n-1}(x)U_{n-2}(y)-U_{n-2}(x)U_{n-1}(y)\big),
\endmultline$$
with $x=(\al+s)/(2\sqrt t)$ and $y=(\be+s)/(2\sqrt t)$.
\endproclaim

\demo{Proof}
The Chebyshev polynomials $U_n(x)$ are defined by
$$
U_j(\cos\th)=\frac {\sin((j+1)\th)} {\sin\th}
=\frac {e^{\bold i(j+1)\th}-e^{-\bold i(j+1)\th}} 
{e^{\bold i\th}-e^{-\bold i\th}}.
\tag\EH
$$
In particular, they satisfy the recurrence
$$
U_{n+1}(x)=2xU_n(x)-U_{n-1}(x), \quad \text{for }n\ge1,
\tag\EHa
$$
with initial values $U_0(x)=1$ and $U_1(x)=2x$.
Hence, the polynomials\linebreak 
$t^{n/2}U_n\big((\al+s)/(2\sqrt t)\big)$
satisfy the recurrence~(\BA) under the given specialisations of the
$s_i$'s and $t_i$'s. Thus, also the right-hand side of (\EF)
satisfies~(\BA). It is then routine to check that it also has
the same initial values as the $f_n(\al)$'s. This proves~(\EF).
The alternative expression (\EFa) (to which we shall come back
in Section~7) follows from this
by an application of~(\EHa).

\medskip
For the second claim, we use the right-most expression in (\EH)
to rewrite $U_j(x)$ as $\frac {z^{j+1}-z^{-j-1}} {z-z^{-1}}$,
where $z=e^{\bold i\th}$ and $x=\cos\th$.
Then the evaluation of the sum on the right-hand side of~(\BC)
amounts to an exercise in computing geometric series and rearranging
terms (which one obviously performs
with the help of a computer algebra programme). 
In the end, one arrives at the right-hand side of~(\EG).\quad \quad \qed
\enddemo

Obviously, we can now take any of the specialisations (i)--(xi)
and (xiv)--(xviii), substitute the corresponding values of
$s_0,s,t_0,t$ in (\EG), and get a Hankel determinant evaluation
in which the entries are linear combinations of three
consecutive members of a well-known sequence. Of course, of special
interest are those cases where the evaluations of the
Chebyshev polynomials in (\EG) have a closed form.
In this regard, we have
$$\align 
U_n(0)&=[n\text{ even}](-1)^{n/2},
\tag\YA\\
U_n(1/2)&=\cases 
1,&\text{if }n\equiv 0,1~\text{mod }6,\\
0,&\text{if }n\equiv 2~\text{mod }3,\\
-1,&\text{if }n\equiv 3,4~\text{mod }6,
\endcases
\tag\YB\\
U_n(1)&=n+1,
\tag\YC\\
U_n(3/2)&=F_{2n+2},
\tag\YD\\
U_n(\bold i/2)&=\bold i^nF_{n+1},
\tag\YE\\
U_n(\bold i)&=\bold i^nP_{n+1},
\tag\YF
\endalign$$
where $P_n$ is the $n$-th {\it Pell number} defined by
$P_n=2P_{n-1}+P_{n-2}$ with initial values $P_0=0$ and $P_1=1$.

\medskip
(I)
In order to obtain (\AA), we have to apply
the specialisation~(ii) in (\CB) with $\al=1$. 
In that case, we have $m_n=C_n$ and
$\frac {\al+s}{2\sqrt t}=\frac {3} {2}$, and therefore
$f_n(1)=F_{2n+2}-F_{2n}=F_{2n+1}$ by (\EF) and (\YD).

In the same vein, to obtain (\AB) we need to apply
the specialisation~(iii) in (\CB) with $\al=1$. 
Here we have $m_n=C_{n+1}$ and
$\frac {\al+s}{2\sqrt t}=\frac {3} {2}$, and therefore
$f_n(\al)=F_{2n+2}$ by (\EF) and (\YD).

Next, in view of the above, we see that application of
specialisation~(ii) in (\BC) with $\al=\be=1$ leads directly to~(\AC),
while specialisation~(iii) produces
$$
\det\left(C_{i+j+1}+2C_{i+j+2}+C_{i+j+3}\right)_{i,j=0}^{n-1}
=\sum_{j=0}^n F_{2j+2}^2.
\tag\YG
$$
On the other hand, if instead we set $\al=\be=-1$ in (\BC) then,
by (\YB), we obtain
$$
\det\left(C_{i+j}-2C_{i+j+1}+C_{i+j+2}\right)_{i,j=0}^{n-1}
=\sum_{j=0}^n \big(U_j(1/2)-U_{j-1}(1/2)\big)^2=
\left\lfloor{\frac {2n+3} {3}}\right\rfloor
\tag\YGa
$$
and
$$
\det\left(C_{i+j+1}-2C_{i+j+2}+C_{i+j+3}\right)_{i,j=0}^{n-1}
=\sum_{j=0}^n U_j^2(1/2)=\left\lfloor{\frac {2n+4} {3}}\right\rfloor.
\tag\YGb
$$

More generally, without specialising $\al$ and $\be$ in (\BC), we
get
$$\multline 
\det\left(\al\be
C_{i+j}+(\al+\be)C_{i+j+1}+C_{i+j+2}\right)_{i,j=0}^{n-1}\\
=\frac {\matrix \big(U_{n+1}(x)U_n(y)-U_n(x)U_{n+1}(y)
+U_{n-1}(x)U_{n+1}(y)\kern3cm\\\kern3cm-U_{n+1}(x)U_{n-1}(y)
+U_{n}(x)U_{n-1}(y)-U_{n-1}(x)U_{n}(y)\big)\endmatrix} {\al-\be},
\endmultline
\tag\AX$$
with $x=(\al+2)/2$ and $y=(\be+2)/2$, and
$$\multline 
\det\left(\al\be
C_{i+j+1}+(\al+\be)C_{i+j+2}+C_{i+j+3}\right)_{i,j=0}^{n-1}\\
=
\frac {U_{n+1}\big((\al+2)/2\big)U_n\big((\be+2)/2\big)
-U_n\big((\al+2)/2\big)U_{n+1}\big((\be+2)/2\big)} {\al-\be}.
\endmultline
\tag\AY$$

\medskip
(II)
The most straightforward specialisation, $s_i=t_i=1$ for all~$i$,
summarised above in Case~(i), when applied in~(\BC), leads to
$$
{\det\left(\al\be M_{i+j}+(\al+\be)M_{i+j+1}+M_{i+j+2}\right)_{i,j=0}^{n-1}} 
=\sum_{j=0}^nU_j\left(\tfrac {\al+1} {2}\right)
U_j\left(\tfrac {\be+1} {2}\right).
\tag\YH
$$
In particular, for $\al=\be=1$, by (\YC) this reduces to
$$
{\det\left(M_{i+j}+2M_{i+j+1}+M_{i+j+2}\right)_{i,j=0}^{n-1}} 
=\sum_{j=0}^n(j+1)^2=\frac {(n+1)(n+2)(2n+3)} {6}.
\tag\YI
$$
On the other hand, if instead we set $\al=\be=-1$, then by (\YA) we get
$$
{\det\left(M_{i+j}-2M_{i+j+1}+M_{i+j+2}\right)_{i,j=0}^{n-1}} 
=\left\lfloor{\frac {n+2} {2}}\right\rfloor.
\tag\YIa
$$
Another interesting specialisation is $\al=\be=2$. In this case,
Identity~(\CB) becomes
$$
{\det\left(2M_{i+j}+M_{i+j+1}\right)_{i,j=0}^{n-1}} 
=F_{2n+2},
\tag\YIb
$$
and (\BC) becomes
$$
{\det\left(4M_{i+j}+4M_{i+j+1}+M_{i+j+2}\right)_{i,j=0}^{n-1}} 
=\sum_{j=0}^nF_{2j+2}^2.
\tag\YIc
$$

(III)
We come to determinants involving central binomial coefficients.
In Case~(iv) above, we have $m_n=\binom {2n}n$ and
$\frac {\al+s}{2\sqrt t}=\frac {\al+2} {2}$, and therefore
$f_0(1)=1$ and
$$f_n(1)=2U_n(3/2)-3U_{n-1}(3/2)=2F_{2n+2}-3F_{2n}=L_{2n},
\quad \text{for }n\ge1$$ 
by (\EF) and (\YD).
Hence, from (\EG) we obtain 
$$\multline
\frac {\det\left(\al\be \binom {2i+2j}{i+j}+
(\al+\be)\binom {2i+2j+2}{i+j+1}+\binom {2i+2j+4}{i+j+2}\right)_{i,j=0}^{n-1}} 
{2^{n-1}}\\
=\frac {\matrix \big(
U_{n+1}(x)U_n(y)-U_n(x)U_{n+1}(y)
+
U_{n}(x)U_{n-1}(y)-U_{n-1}(x)U_{n}(y)
\kern2.5cm\\\kern1cm
-
U_{n+1}(x)U_{n-2}(y)+U_{n-2}(x)U_{n+1}(y)
+
U_{n-1}(x)U_{n-2}(y)-U_{n-2}(x)U_{n-1}(y)
\big)\endmatrix} {\al-\be},
\endmultline
\tag\YJ
$$
with $x=(\al+2)/2$ and $y=(\be+2)/2$. For $\al=\be=1$, 
Identity~(\BC) tells us that the above reduces to
$$
\frac {\det\left( \binom {2i+2j}{i+j}+
2\binom {2i+2j+2}{i+j+1}+\binom {2i+2j+4}{i+j+2}\right)_{i,j=0}^{n-1}} 
{2^{n-1}}
=-2+
\sum_{j=0}^n L_{2j}^2.
\tag\YK
$$
(The ``correction'' of $-2$ on the right-hand side above is caused by the
``discrepancy'' between $f_0(\al)=1$ and $L_0=2$.)

In the same vein, in Case~(xiv) we have 
$m_n=\frac {1} {2}\binom {2n+2}{n+1}$ and
$\frac {\al+s}{2\sqrt t}=\frac {\al+2} {2}$, and therefore
$f_n(1)=U_n(3/2)+U_{n-1}(3/2)=F_{2n+2}+F_{2n}=L_{2n+1}$
by (\EF) and (\YD).
Hence, from (\EG) we obtain 
$$\multline
\frac {
\det\left(\al\be \binom {2i+2j+2}{i+j+1}+
(\al+\be)\binom {2i+2j+4}{i+j+2}+\binom {2i+2j+6}{i+j+3}\right)_{i,j=0}^{n-1}} 
{2^n}\\
=\frac {\matrix \big(
U_{n+1}(x)U_n(y)-U_n(x)U_{n+1}(y)
+
U_{n+1}(x)U_{n-1}(y)\kern3cm\\\kern3cm-U_{n-1}(x)U_{n+1}(y)
+
U_{n}(x)U_{n-1}(y)-U_{n-1}(x)U_{n}(y)\big)\endmatrix
} {\al-\be},
\endmultline
\tag\YL
$$
with $x=(\al+2)/2$ and $y=(\be+2)/2$.
For $\al=\be=1$, 
Identity~(\BC) tells us that the above reduces to~(\AL).

\medskip
(IV)
We continue with central binomial coefficients, but with even {\it
and\/} odd central binomial coefficients, which are obtained by the
choice described in Case~(xv). In this case,
we have $m_n=\binom {n}{\fl{n/2}}$ and
$\frac {\al+s}{2\sqrt t}=\frac {\al} {2}$, and therefore
$$f_n(1)=U_n(1/2)+U_{n-1}(1/2)=\cases 
1,&\text{if }n\equiv0,2~(\text{mod }6),\\
-1,&\text{if }n\equiv3,5~(\text{mod }6),\\
2,&\text{if }n\equiv1~(\text{mod }6),\\
-2,&\text{if }n\equiv4~(\text{mod }6),
\endcases
$$ 
and
$$f_n(-1)=U_n(-1/2)+U_{n-1}(-1/2)=\cases 
1,&\text{if }n\equiv0~(\text{mod }3),\\
-1,&\text{if }n\equiv2~(\text{mod }3),\\
0,&\text{if }n\equiv1~(\text{mod }3),
\endcases
$$ 
by (\EF) and (\YB). Further special evaluations are
$f_n(2)=(n+1)+n=2n+1$ and $f_n(3/2)=F_{2n+2}+F_{2n}=L_{2n+1}$.
Hence, from (\EG) we obtain 
$$\multline
{\det\left(\al\be \binom {i+j}{\fl{(i+j)/2}}+
(\al+\be)\binom {i+j+1}{\fl{(i+j+1)/2}}+
\binom {i+j+2}{\fl{(i+j+2)/2}}\right)_{i,j=0}^{n-1}} 
\\
=\frac {\matrix \big(
U_{n+1}(x)U_n(y)-U_n(x)U_{n+1}(y)
+
U_{n+1}(x)U_{n-1}(y)\kern3cm\\\kern3cm-U_{n-1}(x)U_{n+1}(y)
+
U_{n}(x)U_{n-1}(y)-U_{n-1}(x)U_{n}(y)
\big)\endmatrix} 
{\al-\be},
\endmultline
\tag\YO
$$
with $x=\al/2$ and $y=\be/2$.
In particular, for $\al=\be=1$, this reduces to
$$
\multline
{\det\left(\binom {i+j}{\fl{(i+j)/2}}+
2\binom {i+j+1}{\fl{(i+j+1)/2}}+
\binom {i+j+2}{\fl{(i+j+2)/2}}\right)_{i,j=0}^{n-1}} \\
=\cases 
2n+1,&\text{if }n\equiv0~(\text{mod }3),\\
2n+3,&\text{if }n\equiv1~(\text{mod }3),\\
2n+4,&\text{if }n\equiv2~(\text{mod }3).
\endcases
\endmultline
\tag\YP
$$
On the other hand, for $\al=\be=-1$ we get
$$
{\det\left(\binom {i+j}{\fl{(i+j)/2}}-
2\binom {i+j+1}{\fl{(i+j+1)/2}}+
\binom {i+j+2}{\fl{(i+j+2)/2}}\right)_{i,j=0}^{n-1}} 
=\left\lfloor\frac {2n+3} {3}\right\rfloor.
\tag\YQ
$$
For $\al=2$ and $\al=3$, curious evaluations can be obtained from
(\CB), namely
$$
{\det\left(2\binom {i+j}{\fl{(i+j)/2}}+
\binom {i+j+1}{\fl{(i+j+1)/2}}\right)_{i,j=0}^{n-1}} 
=2n+1
\tag\YR
$$
and
$$
{\det\left(3\binom {i+j}{\fl{(i+j)/2}}+
\binom {i+j+1}{\fl{(i+j+1)/2}}\right)_{i,j=0}^{n-1}} 
=L_{2n+1}.
\tag\YS
$$

\medskip
(V)
Another curious special case arises if we choose $s_i=0$ and $t_i=-1$
for all~$i$. Then $m_n=[n\text{ even}](-1)^{n/2}C_{n/2}$
and $f_n(1)=(-\bold i)^nU_n(\bold i/2)=F_{n+1}$. (The important point
here is to choose $\sqrt{t}=\sqrt{-1}$ {\it consistently}; our choice
here is~$\sqrt{-1}=-\bold i$.)
From (\CB) we therefore get
$$
\multline
\det\left((-1)^{\fl{(i+j+1)/2}}C_{\fl{(i+j+1)/2}}\right)_{i,j=0}^{n-1}=
\det\left(m_{i+j}+m_{i+j+1}\right)_{i,j=0}^{n-1}\\
=
\det\left(m_{i+j}\right)_{i,j=0}^{n-1}\cdot F_{n+1}
=(-1)^{\binom n2}F_{n+1}.
\endmultline
\tag\AU
$$
For example, for $n =5$, we get
$$
\det\pmatrix 
1&-1&-1&2&2\\
-1&-1&2&2&-5\\
-1&2&2&-5&-5\\
2&2&-5&-5&14\\
2&-5&-5&14&14
\endpmatrix
=8=F_6.
$$
Formula (\BC) with the same choices of the $s_i$'s and $t_i$'s gives
$$
\det\left(m_{i+j}+2m_{i+j+1}+m_{i+j+2}\right)_{i,j=0}^{n-1}
=(-1)^{\binom {n+1}2}\sum_{j=0}^n(-1)^jF_{j+1}^2.
\tag\AV
$$

\medskip
(VI)
We point out that our formula (\EG) allows us to vastly generalise 
the second main result in \cite{\DoFSAA, Theorem~3}.
Namely, there Dougherty, French, Saderholm and Qian establish a linear
recurrence of order~4 for the Hankel determinant
$$
{\det\left(\la
C_{i+j}+\mu C_{i+j+1}+C_{i+j+2}\right)_{i,j=0}^{n-1}}  ,
$$
where $\la$ and $\mu$ are some constants. Our formula in (\EG)
provides even a {\it closed form expression} for a determinant
that has much more general matrix entries
(and the reparameterisation $\la=\al\be$ and $\mu=\al+\be$).
Using this, we may deduce the following corollary.

\proclaim{Corollary \TE}
With the setup in Section~2, let $s_i\equiv s$ and $t_i\equiv t$ for
$i\ge1$. Then the scaled Hankel determinants 
$$
t_0^{-(n-1)}t^{-\binom {n-1}2}
\det\left(\al\be m_{i+j}+(\al+\be)m_{i+j+1}+m_{i+j+2}\right)_{i,j=0}^{n-1}
\tag\YM
$$
satisfy a linear recurrence of order~$4$ with constant coefficients.
\endproclaim

\remark{Remark}
The reader should note that the scaling in (\YM) is exactly the
value of the determinant in the denominator on the left-hand side of (\EG). 
\endremark

\demo{Proof of Corollary \TE}
The reader should note that the expression on the right-hand side
of~(\EG) is a 
linear combination of products of the form
$t^nU_{n+\de}(x)U_{n+\ep}(y)$ with constant coefficients,
where $\de,\ep\in\{1,0,-1,-2\}$. (We recall that the arguments of the
Chebyshev polynomials are 
$x=(\al+s)/(2\sqrt t)$ and $y=(\be+s)/(2\sqrt t)$.) Now, by (\EHa),
for each fixed~$\de$ the sequence $\big(t^{n/2}U_{n+\de}(x)\big)_{n\ge0}$
satisfies the recurrence relation
$$
p_{n+1}=(\al+s)p_n-tp_{n-1}.
\tag\YN
$$
Similarly, for each fixed~$\ep$ the sequence 
$\big(t^{n/2}U_{n+\ep}(y)\big)_{n\ge0}$
satisfies the same recurrence relation with $\al$ replaced by~$\be$.
Hence, each product sequence $\big(t^nU_{n+\de}(x)U_{n+\ep}(y)\big)_{n\ge0}$
satisfies the same recurrence relation,
namely the one resulting from the ``product''
of (\YN) with (\YN) where $\al$ has been replaced by~$\be$.
It follows from the proof of \cite{\StanBI, Theorem~6.4.9}
(which is actually a much more general theorem) that the order of
this ``product'' recurrence is $2\cdot2=4$.\quad \quad \qed
\enddemo

\medskip
(VII) We turn to $q$-analogues. We use the familiar notations
$(a;q)_0:=1$, $(a;q)_n:=(1-a)(1-aq)\cdots(1-aq^{n-1})$, $n\ge1$,
for the {\it $q$-shifted factorials}, and
$$
\bmatrix n\\k\endbmatrix_q:=\frac {(q;q)_n} {(q;q)_k\,(q;q)_{n-k}}
$$
for the {\it $q$-binomial coefficients}.
The {\it Rogers--Szeg\H o polynomials}
$r_n(t)$, defined by
$$
r_n(t):=\sum_{k=0}^n\bmatrix n\\k\endbmatrix_qt^k,
$$
are the moments for the orthogonal polynomials (\BBa) with
$s_i=q^i(t+1)$ and $t_i=\mathbreak q^it(q^{i+1}-1)$. The corresponding
polynomials $f_n(\al)$ are given by
$$
f_n(\al)=\sum_{k=0}^n\bmatrix n\\k\endbmatrix_q
\sum_{j=0}^kq^{\binom j2+\binom {k-j}2}\bmatrix k\\j\endbmatrix_q
t^j\al^{n-k}.
$$
Thus, from (\CB) we obtain
$$
\frac {\det\left(\al r_{i+j}(t)+r_{i+j+1}(t)\right)_{i,j=0}^{n-1}} 
{(-t)^{\binom n2}q^{\binom n3}
\prod _{i=1} ^{n-1}(q,q)_i}
=\sum_{k=0}^n\bmatrix n\\k\endbmatrix_q
\sum_{j=0}^kq^{\binom j2+\binom {k-j}2}\bmatrix k\\j\endbmatrix_q
t^j\al^{n-k},
\tag\AK$$
and from (\BC) we get
$$\multline 
\frac {\det\left(\al\be r_{i+j}(t)+(\al+\be)t_{i+j+1}(t)+r_{i+j+2}(t)\right)_{i,j=0}^{n-1}} 
{(-t)^{\binom n2}q^{\binom n3}
\prod _{i=1} ^{n-1}(q,q)_i}\\
=
\sum_{j=0}^nf_{j}(\al)f_j(\be)
(-t)^{n-j}q^{\binom n2-\binom j2}(q^{j+1};q)_{n-j},
\endmultline
\tag\AM$$
with $f_j(\al)$ as above. We point out that the right-hand sides in
(\AK) and (\AM) simplify for $t=-1$ since, in basic hypergeometric
notation (cf\. \cite{\GaRaAF}), we have
$$
\sum_{j=0}^kq^{\binom j2+\binom {k-j}2}\bmatrix k\\j\endbmatrix_q
t^j=
q^{\binom k2}
  {} _{1} \phi _{1} \! \left [ \matrix q^{-k}\\ 0\endmatrix ;q, qt \right ]  ,
$$
and this $_1\phi_1$-series can be evaluated by means of \cite{\GaRaAF,
(1.8.1); Appendix (II.9); $b\to\infty$}. To be precise, for $t=-1$ we
have
$$
\sum_{j=0}^kq^{\binom j2+\binom {k-j}2}\bmatrix k\\j\endbmatrix_q
(-1)^j=\cases 
(-1)^{k/2}q^{\frac {1} {4}(k^2-2k)}(q;q^2)_{k/2},
&\text{if $k$ is even,}\\
0,&\text{if $k$ is odd,}
\endcases
$$
and thus for example
$$
\frac {\det\left(\al r_{i+j}(-1)+r_{i+j+1}(-1)\right)_{i,j=0}^{n-1}} 
{q^{\binom n3}
\prod _{i=1} ^{n-1}(q,q)_i}
=\sum_{k=0}^n\bmatrix n\\2k\endbmatrix_q
(-1)^{k}q^{k^2-k}(q;q^2)_{k}\al^{n-2k},
\tag\AMa
$$
by (\AK).
For $\al=1$, this reduces even further, namely to
$$
\frac {\det\left(r_{i+j}(-1)+r_{i+j+1}(-1)\right)_{i,j=0}^{n-1}} 
{q^{\binom n3}
\prod _{i=1} ^{n-1}(q,q)_i}
=q^{\binom n2},
\tag\AMb
$$
which is most easily seen by looking at the recursion (\BA)
for $\al=1$, $s_i=0$, and $t_i=t^i(t^{i+1}-1)$ for all~$i$.

\medskip
(VIII) The powers $q^{\binom n2}$ are 
the moments for the orthogonal polynomials (\BBa) with
$s_i=q^{i-1}(q^{i+1}+q^i-1)$ and $t_i=q^{3i}(q^{i+1}-1)$. The corresponding 
polynomials $f_n(\al)$ are given by
$$
f_n(\al)=\sum_{k=0}^n\bmatrix n\\k\endbmatrix_q
q^{(n-1)k}\al^{n-k}.
$$
Thus, from (\BC) we obtain
$$
\multline
\frac {\det\left(\al\be q^{\binom {i+j}2}+(\al+\be)q^{\binom {i+j+1}2}
+q^{\binom {i+j+2}2}\right)_{i,j=0}^{n-1}} 
{(-1)^{\binom n2}q^{3\binom n3}
\prod _{i=1} ^{n-1}(q;q)_i}\\
=\sum_{j=0}^n
(-1)^{n-j}q^{3\binom n2-3\binom j2}(q^{j+1};q)_{n-j}
\sum_{k_1=0}^j\bmatrix j\\k_1\endbmatrix_q q^{(j-1)k_1}\al^{j-k_1}
\sum_{k_2=0}^j\bmatrix j\\k_2\endbmatrix_q q^{(j-1)k_2}\be^{j-k_2}.
\endmultline
\tag\AO
$$

\medskip
(IX) For another example, we choose $T_i=q^{\fl{(i+1)/2}}
\left(1-q^{\fl{(i+2)/2}}\right)$ in the Dyck path setting in
Section~2. Then we get
$c_{n}=(q;q)_{n/2}$ for even~$n$ and $c_n=0$ for odd~$n$. 
The corresponding polynomials are
$$
g_{2n}(\al)=\sum_{k=0}^n q^{\binom k2}(q;q)_k\bmatrix
n\\k\endbmatrix_q^2\al^{n-k}
$$
and
$$
g_{2n+1}(\al)=\sum_{k=0}^n q^{\binom k2}(q;q)_k
{\bmatrix n+1\\k\endbmatrix_q}
{\bmatrix n\\k\endbmatrix_q}
\al^{n-k}.
$$
Hence, from (\EC) with $\al=1$ we obtain
$$
\frac {\det\left((q;q)_{\fl{(i+j+1)/2}}\right)_{i,j=0}^{n-1}} 
{q^{\frac {1} {6} n (2 n^2-3n+1)}
\prod _{i=1} ^{n-1}(q;q)_i^2}
=\sum_{k=0}^n q^{\binom k2}(q;q)_k\bmatrix
n\\k\endbmatrix_q^2,
\tag\AR
$$
and similarly, from (\ED),
$$
\frac {\det\left((q;q)_{\fl{(i+j+2)/2}}\right)_{i,j=0}^{n-1}} 
{q^{2\binom {n+1}3}(q;q)_n
\prod _{i=1} ^{n-1}(q;q)_i^2}
=\sum_{k=0}^n q^{\binom k2}(q;q)_k\bmatrix n+1\\k\endbmatrix_q
\bmatrix n\\k\endbmatrix_q.
\tag\AS
$$

\medskip
(X) Obviously, our identities in Sections~2
and~3 lead to many more Hankel determinant evaluations 
involving well-known functions
when applied to the {\it ``orthogonal polynomials as moments''}
from \cite{\IsStAB, \IsStAC}.

\subhead 5. Computational proofs using matrix algebra\endsubhead

\demo{Proof of Theorem \UA}
Define the matrices $G_n(\al)$ by
$$
G_n(\al):=\pmatrix 
\al+s_0&1& \\
t_0&\al+s_1&1&\\
&t_1&\al+s_2&1&\\
\hdotsfor5 \\
&&t_{n-3}&\al+s_{n-2}&1\\
&&&t_{n-2}&\al+s_{n-1}
\endpmatrix
$$
(we show only non-zero entries in the matrix).
By expanding the determinant along the last row of the matrix,
we obtain
$$
\det G_n(\al)
=(\al+s_{n-1})\det G_{n-1}(\al)-t_{n-2}\det G_{n-2}(\al).
$$
Since we have in addition 
$$\det G_1(\al)=\al+s_0=f_1(\al)\quad \text{and}\quad 
\det G_2(\al)=(\al+s_0)(\al+s_1)-t_0=f_2(\al),$$
we see that
$$
f_n(\al)=\det G_n(\al).
\tag\AN
$$

Let $A_n$ be the matrix of path generating functions
$\big(m(i,j)\big)_{i,j=0}^{n-1}$, 
and let $D_n$ be the diagonal matrix with entries
$d(i,i)=
\prod _{\ell=0} ^{i-1}t_\ell$.
By (\BB) we get
$$
\big(m(i+1,j)\big)_{i,j=0}^{n-1}=A_nR_n,
\tag\CC$$
where
$$
R_n:=
\pmatrix 
s_0&1& \\
t_0&s_1&1&\\
&t_1&s_2&1&\\
\hdotsfor5 \\
&&t_{n-3}&s_{n-2}&1\\
&&&t_{n-2}&s_{n-1}
\endpmatrix.$$
It is well known (cf\. e.g\. \cite{\CiglAT}) and easy to see by
cutting a Motzkin path from $(0,0)$ to $(i+j,0)$ after $i$ steps that
$$
\sum_{k}m(i,k)m(j,k)
\prod _{\ell=0} ^{k-1}t_\ell=m(i+j,0).
\tag\CD
$$
This means that
$$\big(m_{i+j}\big)_{i,j=0}^{n-1}=A_nD_nA_n^t,$$ 
and it implies that
$$
\big(\al m_{i+j}+m_{i+j+1}\big)_{i,j=0}^{n-1}
=\al A_nD_nA_n^t+A_nR_nD_nA_n^t=A_nG_n(\al)D_nA_n^t.
$$
Since $A_n$ is a lower diagonal matrix and 
$m(i,i)=1$ for all $i$, we get 
$$\det\big( A_nG_n(\al)D_nA_n^t\big)=\det D_n\cdot\det G_n(\al),$$
which, recalling (\BD) and (\AN), implies (\CB).

In the same way as before (\CC), we get
$$
\big(m(i+2,j)\big)_{i,j=0}^{n-1}=A_n(R_n^2+t_{n-1}E_n),
\tag\CE$$
where $E_n=\big(u_n(i,j)\big)_{i,j=0}^{n-1}$ with $u_n(n-1,n-1)=1$ and
$u_n(i,j)=0$ otherwise. The term $t_{n-1}E_n$ arises to account for
the Motzkin paths that are included in the weighted enumeration 
$m(n+1,n-1)$ which end with an up-step
followed by a down-step.
Then, using (\CD) with $i$ replaced by $i+2$, we infer
$$
\align
\big(\al \be m(i+j,0)+(\al+\be)&m(i+j+1,0)+m(i+j+2,0)\big)_{i,j=0}^{n-1}\\
&=A_n(\al\be I_n+(\al+\be)R_n+R_n^2+t_{n-1}E_n)D_nA_n^t\\
&=A_n(G_n(\al)G_n(\be)+t_{n-1}E_n)D_nA_n^t,
\tag\CF
\endalign
$$
where $I_n$ denotes the $n\times n$ identity matrix.

If\/ $[A]_k$ denotes the restriction of the matrix $A$ to the first $k$ rows
and columns, then it is a simple exercise to verify that
$$
\big[G_n(\al)G_n(\be)\big]_{n-1}
=G_{n-1}(\al)G_{n-1}(\be)+t_{n-2}E_{n-1}.
$$
If we write $H_n$ for the matrix $G_n(\al)G_n(\be)+t_{n-1}E_n$, then
the above relation entails a recursion for the determinant of
this matrix, namely
$$\det H_n=\det\big( G_n(\al)G_n(\be)\big)
+t_{n-1}\det [G_n(\al)G_n(\be)]_{n-1}
=f_n(\al)f_n(\be)+t_{n-1}\det H_{n-1}.
$$
Consequently,
by induction, we get $\det H_n=\sum_{j=0}^nf_j(\al)f_j(\be)
\prod _{\ell=j} ^{n-1}t_\ell$.
If this is substituted in (\CF) and the determinant is taken on both
sides of that identity, then we obtain (\BC) if we again
recall (\BD).\quad \quad \qed
\enddemo

\demo{Proof of Corollary \TC}
There is a well-known correspondence between Dyck paths and Motzkin
paths which arises by arranging the steps of a Dyck path in groups of two,
and then interpreting a group of two up-steps as an up-step in the
corresponding Motzkin path, a group of two down-steps as a down-step in the
corresponding Motzkin path, and a group consisting of an up-step and a
down-step as a level step in the corresponding Motzkin path.
This idea stands behind all the computations below.

With the $c(n,k)$'s defined by (\AE), we let
$a(n, k) =c(2n, 2k)$ and $b(n, k) =c(2n +1, 2k +1)$.
Then, by iterating (\AE), we get
$$
a(n,k)=a(n-1,k-1)+s_ka(n-1,k)+t_ka(n-1,k+1)
\tag\CG
$$
with $s_k=T_{2k-1}+T_{2k}$ and $t_k=T_{2k}T_{2k+1}$, where, by
convention, $T_{-1}=0$,
and
$$
b(n,k)=b(n-1,k-1)+s_kb(n-1,k)+t_kb(n-1,k+1)
\tag\CH
$$
with $s_k=T_{2k}+T_{2k+1}$ and $t_k=T_{2k+1}T_{2k+2}$.

Then, with the $g_n(\al)$'s as defined by (\AD), we get
$$\align
g_{2n}(\al)&=\al g_{2n-1}(\al)+T_{2n-2}g_{2n-2}(\al)\\
&=\al g_{2n-2}(\al)+T_{2n-2}g_{2n-2}(\al)+\al T_{2n-3}g_{2n-3}(\al)\\
&=\al g_{2n-2}(\al)+T_{2n-2}g_{2n-2}(\al)+T_{2n-3}(g_{2n-2}(\al)-T_{2n-4}g_{2n-4}(\al))\\
&=(\al+T_{2n-2}+T_{2n-3})g_{2n-2}(\al)-T_{2n-3}T_{2n-4}g_{2n-4}(\al).
\endalign$$
This implies that for $a(n, k)$ the associated polynomials $f_n(\al)$ 
(in the sense of the setup in Section~2, where the $s_k$'s and $t_k$'s
are as given below (\CG))
are given by $f_n(\al)=g_{2n}(\al)$.

Similarly, we have
$$\align
g_{2n+1}(\al)&=g_{2n}(\al)+T_{2n-1}g_{2n-1}(\al)\\
&=\al g_{2n-1}(\al)+T_{2n-1}g_{2n-1}(\al)+T_{2n-2}g_{2n-2}(\al)\\
&=\al g_{2n-1}(\al)+T_{2n-1}g_{2n-1}(\al)+T_{2n-2}(g_{2n-1}(\al)-T_{2n-3}g_{2n-3}(\al))\\
&=(\al+T_{2n-1}+T_{2n-2})g_{2n-1}(\al)-T_{2n-2}T_{2n-3}g_{2n-3}(\al).
\endalign$$
This implies that for $b(n, k)$  the associated polynomials $f_n(\al)$ 
(defined by the $s_k$'s and $t_k$'s as given below (\CH))
are given by $f_n(\al)=g_{2n+1}(\al)$.

Therefore, Identity~(\BC) with $s_k=T_{2k-1}+T_{2k}$ and
$t_k=T_{2k}T_{2k+1}$ directly gives~(\AI) since in that
case we have $m_n=a(n,0)=c(2n,0)=c_{2n}$. 
Similarly,
by observing that $c(2n+2,0) =T_0c(2n +1,1)$ for $n\ge0$, we see that
Identity~(\BC) with $s_k=T_{2k}+T_{2k+1}$ and $t_k=T_{2k+1}T_{2k+2}$ 
gives~(\AJ) since in that 
case we have $m_n=b(n,0)=c(2n+1,1)=T_0^{-1}c(2n+2,0)=
T_0^{-1}c_{2n+2}$.\quad \quad \qed
\enddemo

\subhead 6. Combinatorial proofs using non-intersecting lattice paths
\endsubhead

\demo{Proof of Theorem \UA}
In order to give a combinatorial proof of
(\BC), we need a combinatorial model for
the polynomials $f_n(\al)$, and we need (another) combinatorial model
for the Hankel determinant 
$$\det\left(\al\be m_{i+j}+(\al+\be)m_{i+j+1}+m_{i+j+2}
\right)_{i,j=0}^{n-1}.
\tag\DA
$$

Our model for $f_n(\al)$ consists of $n\times 1$ ``pillars''
which are decomposed into ``squares''
($1\times 1$ bricks) --- which come in two kinds --- 
and ``dominoes'' ($2\times 1$ bricks).
An example for $n=8$ can be found in Figure~\FA.
There, squares of the second kind are indicated by thick
framing, whereas squares of the first kind are shown ``normally".
(The labelling should be ignored at this point.)
This pillar consists --- from bottom to top --- of a square of
the second kind, followed
by two dominoes, followed by a square of the first kind, 
a square of the second kind, and a square of the first kind.

\midinsert
$$
\Einheit.6cm
\PfadDicke{1pt}
\Pfad(0,0),122222222\endPfad
\Pfad(0,0),222222221\endPfad
\Pfad(0,1),1\endPfad
\Pfad(0,3),1\endPfad
\Pfad(0,5),1\endPfad
\Pfad(0,6),1\endPfad
\Pfad(0,7),1\endPfad
\PfadDicke{3pt}
\Pfad(0,0),1256\endPfad
\Pfad(0,6),1256\endPfad
\Label\ro{\al}(0,0)
\Label\ro{-t_1}(0,1)
\Label\ro{-t_3}(0,3)
\Label\ro{s_5}(0,5)
\Label\ro{\al}(0,6)
\Label\ro{s_7}(0,7)
\hskip.6cm
$$
\centerline{\eightpoint Figure \FA}
\endinsert

We define the weight of a square of the first kind 
whose bottom is at height~$h$
to be $s_h$, we define the weight of a square of the second kind 
to be $\al$,
and we define the weight of a domino whose bottom
is at height~$h$ to be $-t_h$. The weight of a pillar $p$ is by definition
the product of the weights of its squares and dominoes and is denoted
by $w(p)$.
Thus, the weight of the pillar in Figure~\FA\ is
$$
\al(-t_1)(-t_3)s_5\al s_7
$$
(which now also explains the labels in Figure~\FA).

Let $\Cal P_n$ denote the set of all $n\times 1$ pillars.
We claim that
$$
\sum_{p\in \Cal P_n}w(p)=f_n(\al).
\tag\DD
$$
This is indeed easy to see. For, on top of such a pillar we may
either have a square of either kind,
where the sum of the weights of these two possible squares 
is $\al+s_{n-1}$, or a
domino, which has weight $-t_{n-2}$. If we remove the top brick,
then there remains an $(n-1)\times1$, respectively an $(n-2)\times1$
pillar. Consequently, the generating function for $n\times1$ pillars,
$\sum_{p\in \Cal P_n}w(p)$, satisfies as well the recurrence
(\BA), and the initial conditions are obvious. This establishes
the claim.

\medskip
We turn to our model for the Hankel determinant (\DA).
Not very surprisingly, it is a model of non-intersecting lattice paths
and is based on the corresponding enumeration theorem due to
Lindstr\"om \cite{\LindAA, Lemma~1} 
and Gessel and Viennot \cite{\GeViAA; \GeViAB, Theorem~1}.
We recall that a family of paths in a (directed) graph is called
 {\it non-intersecting} if no two paths of the family
have a {\it graph vertex} in common.

The underlying directed graph is the following: vertices are the
points $(x,y)$ in the two-dimensional integer lattice $\Z\times \Z$
with non-negative ordinate (i.e., with $y\ge0$). The
edges are {\it up-steps} $(x,y)\to (x +1,y +1)$ with weight 1, {\it down-steps}
$(x,y)\to (x +1,y -1)$
with weight $t_{y-1}$ if $y \ge1$, and {\it horizontal steps} $(x,y)\to (x
+1,y)$ with
weight $s_y$; in a strip of width~2 centred at the $y$-axis,
we introduce further horizontal steps of a ``second kind'':
in our directed graph
we have {\it additional\/} 
horizontal steps $(-1,y)\to (0,y)$ with weight $\al$,
and {\it additional\/} horizontal steps $(0,y)\to (1,y)$ with weight $\be$.
As before, 
the weight $w(P)$ of a path $P$ is the product of the weights of all
its edges. 

Let $e(A_i,E_j):=\sum_{P:A_i\to E_j}w(P)$, 
where the sum is over all paths from $A_i=(-i-1,0)$ to $E_j=(j+1,0)$,
for $i,j\ge 0$.
Then we have
$$
e(A_i,E_j)=\al\be m_{i+j}+(\al+\be)m_{i+j+1}+m_{i+j+2}.
$$
For, the set of those paths which do not contain any of the additional
horizontal steps contributes $m_{i+j+2}$ to the total weight,
the set of those paths which use one of the additional horizontal
steps of the type $(-1,y)\to (0,y)$ contributes $\al m_{i+j+1}$, the set of
those paths which use one of the additional horizontal
steps of the type $(0,y)\to (1,y)$ contributes $\be m_{i+j+1}$, and the
set of those paths which which use both types of additional horizontal
steps contributes $\al\be m_{i+j}$.

The Lindstr\"om--Gessel--Viennot theorem then gives
$$
\multline
\det\left(\al\be m_{i+j}+(\al+\be)m_{i+j+1}+m_{i+j+2}\right)_{i,j=0}^{n-1}\\
=
\det\left(e(A_i,E_j)\right)_{i,j=0}^{n-1}
=
\sum_{\Cal P}\sgn(\si(\Cal P))
\prod _{i=0} ^{n-1}w(P_i).
\endmultline
\tag\AW
$$
Here, $\Cal P=(P_0,P_1,\dots,P_{n-1})$ runs over all families of
non-intersecting
paths, where $P_i$ runs from $A_i$ to $E_{\si(i)}$, $i=0,1,\dots,n-1$, 
for some permutation $\si$ of $\{0,1,\dots,n-1\}$.
An example of such a family of non-intersecting lattice paths
for $n=9$ is shown in Figure~\FB. There, the additional horizontal
steps between the abscissa $-1$ and $+1$ are indicated by thick
segments. For example, path $P_0$ consists of two of these additional
steps (hence its weight is $\al\be$), while path $P_6$ contains
an additional horizontal step with weight $\al$, and path $P_7$ contains
an additional horizontal step with weight $\be$.
The permutation $\si$ associated with the path family in Figure~\FB\
is $031426578$ (meaning that $P_0$ ends in $E_0$, $P_1$ ends in $E_3$,
etc.). 

\midinsert
$$
\Einheit.6cm
\Gitter(11,11)(-10,-1)
\Koordinatenachsen(11,11)(-10,-1)
\SPfad(-1,-2),2222222222222\endSPfad
\SPfad(1,-2),2222222222222\endSPfad
\Pfad(-9,0),333333333444444444\endPfad
\Pfad(-8,0),3333333114444444\endPfad
\Pfad(-7,0),3333331444444\endPfad
\Pfad(-6,0),3333313444444\endPfad
\Pfad(-5,0),33334444\endPfad
\Pfad(-4,0),333314444\endPfad
\Pfad(-3,0),33414\endPfad
\Pfad(-2,0),333444\endPfad
\Pfad(-1,0),11\endPfad
\PfadDicke{3pt}
\Pfad(-1,0),11\endPfad
\Pfad(-1,6),1\endPfad
\Pfad(0,7),1\endPfad
\DickPunkt(-9,0)
\DickPunkt(-8,0)
\DickPunkt(-7,0)
\DickPunkt(-6,0)
\DickPunkt(-5,0)
\DickPunkt(-4,0)
\DickPunkt(-3,0)
\DickPunkt(-2,0)
\DickPunkt(-1,0)
\DickPunkt(1,0)
\DickPunkt(2,0)
\DickPunkt(3,0)
\DickPunkt(4,0)
\DickPunkt(5,0)
\DickPunkt(6,0)
\DickPunkt(7,0)
\DickPunkt(8,0)
\DickPunkt(9,0)
\Label\u{A_0}(-1,0)
\Label\u{A_1}(-2,0)
\Label\u{A_2}(-3,0)
\Label\u{A_3}(-4,0)
\Label\u{A_4}(-5,0)
\Label\u{A_5}(-6,0)
\Label\u{A_6}(-7,0)
\Label\u{A_7}(-8,0)
\Label\u{A_8}(-9,0)
\Label\u{E_0}(1,0)
\Label\u{E_1}(2,0)
\Label\u{E_2}(3,0)
\Label\u{E_3}(4,0)
\Label\u{E_4}(5,0)
\Label\u{E_5}(6,0)
\Label\u{E_6}(7,0)
\Label\u{E_7}(8,0)
\Label\u{E_8}(9,0)
$$
\centerline{\eightpoint Figure \FB}
\endinsert

A reader who is experienced in the combinatorics of non-intersecting
lattice paths does not have to go through the following explanations
in detail; it will suffice to look at Figures~\FB\ and~\FC, since
they contain the essence of the argument. 
Everything else is --- so-to-speak --- just ``book-keeping''.

However, in order to properly explain, we start with the question of
how such an above family of non-intersecting paths looks like?

In the figure, we have marked the vertical lines $x=-1$ and $x=+1$
by dotted lines. To the left of $x=-1$, the paths are uniquely
determined by the condition of being non-intersecting, and the
same is true to the right of $x=+1$; see the figure.
Thus, it is only in the strip between the vertical lines $x=-1$ 
and $x=+1$ where something ``interesting'' may happen.

In that strip, there may be some paths on the top, say $P_{n-1},\dots,P_j$,
which, in that
strip, have an up-step followed by a down-step. (In the figure,
this is only the case for the top-most path $P_8$.) 
All other paths do not exceed height~$j-1$, and they form patterns
in that strip that can be described as follows.
We first look into the left substrip between the vertical lines
$x=-1$ and $x=0$. There, a path $P_i$ may do three different things
(see Figure~\FB):

\roster 
\item a horizontal step, which may be of two different kinds ---
``ordinary'' or one of the ``additional'' steps;
\item a down-step, but then there must be another path
doing an up-step in the substrip crossing $P_i$;
\item an up-step, but then there must be another path
doing a down-step in the substrip crossing $P_i$.
\endroster

The same is true in the substrip between the vertical lines $x=0$
and $x=+1$, and what happens there is {\it independent} from
what is happening in the other substrip.

In order to determine the contribution of all this to the 
overall generating function (\AW), we observe that the weight
contribution of the path parts to the left of the vertical line
$x=-1$ and to the right of the vertical line $x=+1$ is in total
the right-hand side of (\BD), or, in other words,
$$
{\det\left(m_{i+j}\right)_{i,j=0}^{n-1}},
$$
the denominator in (\BC). Hence, the contribution of the various
configurations between the vertical lines $x=-1$ and $x=+1$
must give the right-hand side of (\BC).

This contribution of the ``free'' part in the strip between the vertical lines
$x=-1$ and $x=+1$ is composed of the following constituents:

\roster 
\item The top-most paths $P_{n-1},\dots,P_j$, each of which have 
an up-step followed by a down-step in this strip, contribute
a multiplicative factor of
$$
t_{j}t_{j+1}\cdots t_{n-1}
\tag\DB
$$
to the total weight.
\item According to the above considerations,
the contribution of the paths $P_0,P_1,\dots,\mathbreak P_{j-1}$ below
equals the product of the their contributions in the 
left substrip between the vertical lines $x=-1$ and $x=0$
and their contributions in the
right substrip between the vertical lines $x=0$ and $x=+1$.
Combinatorially, we find the same configurations in the left substrip
and the right substrip. The only difference is that the ``additional''
horizontal steps have weight $\al$ in the left substrip, whereas
they have weight $\be$ in the right substrip.
Hence, if we denote the generating function of the configurations
in the left substrip by $\GF_j(\al)$, then the total contribution
coming from the strip between the vertical lines $x=-1$ and $x=+1$
equals 
$$
\GF_j(\al)\GF_j(\be).
\tag\DC$$
\endroster

Let us concentrate on the contributions of $P_0,P_1,\dots,P_{j-1}$
in the left substrip, that is on $\GF_j(\al)$. In order to see
that $\GF_j(\al)=f_j(\al)$, we set up a weight-preserving bijection
between these contributions and the pillars that we introduced
at the beginning of this section. In this bijection, we map
an ``ordinary'' horizontal step on height~$h$ to a square
of the first kind on height~$h$, we map 
an ``additional'' horizontal step on height~$h$ to a square
of the second kind on height~$h$, and we map a ``cross'', that is
a pair consisting of a down-step and an up-step that cross each
other, with the down-step ending on height~$h$, to a domino on
height~$h$. The left part of Figure~\FC\ shows this correspondence for the
configuration in the left substrip of Figure~\FB.
(For the sake of completeness, the right part
of Figure~\FC\ shows the correspondence for the configuration in
the right substrip of Figure~\FB. The reader is alerted that
the weights of the squares of the second kind must be defined here to
equal~$\be$ in order to set up a weight-preserving
correspondence.)

\midinsert
$$
\Einheit .6cm
\Gitter(1,8)(-1,0)
\Koordinatenachsen(1,8)(-1,0)
\SPfad(-1,-1),222222222\endSPfad
\Pfad(-1,7),1\endPfad
\Pfad(-1,6),1\endPfad
\Pfad(-1,5),1\endPfad
\Pfad(-1,4),4\endPfad
\Pfad(-1,3),3\endPfad
\Pfad(-1,2),4\endPfad
\Pfad(-1,1),3\endPfad
\Pfad(-1,0),1\endPfad
\DickPunkt(-1,0)
\Label\ro{\longleftrightarrow}(1,5)
\PfadDicke{1pt}
\Pfad(3,0),122222222\endPfad
\Pfad(3,0),222222221\endPfad
\Pfad(3,1),1\endPfad
\Pfad(3,3),1\endPfad
\Pfad(3,5),1\endPfad
\Pfad(3,6),1\endPfad
\Pfad(3,7),1\endPfad
\PfadDicke{3pt}
\Pfad(-1,0),1\endPfad
\Pfad(-1,6),1\endPfad
\Label\ro{\al}(3,0)
\Label\ro{-t_1}(3,1)
\Label\ro{-t_3}(3,3)
\Label\ro{s_5}(3,5)
\Label\ro{\al}(3,6)
\Label\ro{s_7}(3,7)
\hbox{\hskip5cm}
\Gitter(2,8)(0,0)
\Koordinatenachsen(2,8)(0,0)
\SPfad(1,-1),222222222\endSPfad
\Pfad(0,7),1\endPfad
\Pfad(0,6),4\endPfad
\Pfad(0,5),3\endPfad
\Pfad(0,4),1\endPfad
\Pfad(0,3),4\endPfad
\Pfad(0,2),3\endPfad
\Pfad(0,1),1\endPfad
\Pfad(0,0),1\endPfad
\DickPunkt(1,0)
\Label\ro{\longleftrightarrow}(2,5)
\PfadDicke{1pt}
\Pfad(4,0),122222222\endPfad
\Pfad(4,0),222222221\endPfad
\Pfad(4,1),1\endPfad
\Pfad(4,2),1\endPfad
\Pfad(4,4),1\endPfad
\Pfad(4,5),1\endPfad
\Pfad(4,7),1\endPfad
\PfadDicke{3pt}
\Pfad(0,0),1\endPfad
\Pfad(0,7),1\endPfad
\Label\ro{\be}(4,0)
\Label\ro{s_1}(4,1)
\Label\ro{-t_2}(4,2)
\Label\ro{s_4}(4,4)
\Label\ro{-t_5}(4,5)
\Label\ro{\be}(4,7)
\hskip2cm
$$
\centerline{\eightpoint Figure \FC}
\endinsert

This is indeed a weight-preserving bijection. For the correspondence
between horizontal steps and squares this is
obvious, but it also works for the correspondence between
crosses and dominoes: in a cross two paths change sides\footnote{It is
important to notice that a cross is not a violation of the property
of being non-intersecting: the ``intersection'' happens in a 
{\it non-lattice point}, that is, the corresponding paths do {\it
not\/} have a vertex of the underlying graph in common.}; in
the corresponding permutation $\si(\Cal P)$ this corresponds to a
transposition, which has negative sign. Indeed, we have defined
the weight of a domino on height~$h$ by $-t_h$, that is, with
a negative sign in front. If we now remember (\DD), then we
see that the proof of (\BC) is complete.\quad \quad \qed
\enddemo

\demo{Sketch of proof of Corollary \TC}
The two Hankel determinant evaluations in Corollary~\TC\ can also
be directly established by a non-intersecting lattice path model.
Since the corresponding arguments are very similar to the previous
proof, we content ourselves with giving just a rough sketch that
essentially only consists in presenting the figures analogous to
the ones in Figures~\FB\ and~\FC. 

Here, the underlying graph is again the upper two-dimensional integer
lattice. However, in difference to the previous proof, there are no
``ordinary'' horizontal steps at all, while the ``additional'' 
horizontal steps are horizontal steps of length~2 from $(-2,y)$ to
$(0,y)$ and from $(0,y)$ to $(2,y)$, for $y\ge0$. The former
horizontal steps have weight~$\al$, while the latter have
weight~$\be$.
As for the ``ordinary'' steps, up-steps have weight~1, and
down-steps $(x,y+1)\to(x+1,y)$ have weight~$T_y$.

\medskip
The non-intersecting lattice path model for the Hankel determinant in
the numerator of (\AI) consists of families
$\Cal P=(P_0,P_1,\dots,P_{n-1})$ of non-intersecting lattice paths
where $P_i$ runs from $A_i=(-2i-2)$ to $E_{\si(i)}$, with
$E_i=(2i+2,0)$ and $\si$ being some permutation of
$\{0,1,\dots,n-1\}$.
See Figure~\FD\ for an example.

\midinsert
$$
\Gitter(11,11)(-10,-1)
\Koordinatenachsen(11,11)(-10,-1)
\SPfad(-2,-2),2222222222222\endSPfad
\SPfad(2,-2),2222222222222\endSPfad
%\Pfad(-14,0),3333333333333344444444444444\endPfad
%\Pfad(-12,0),333333333333444444444444\endPfad
\Pfad(-10,0),33333333334444444444\endPfad
\Pfad(-8,0),3333331143444444\endPfad
\Pfad(-6,0),333343114444\endPfad
\Pfad(-4,0),33433444\endPfad
\PfadDicke{3pt}
\Pfad(-2,0),1111\endPfad
\Pfad(-2,6),11\endPfad
\Pfad(0,4),11\endPfad
%\DickPunkt(-14,0)
%\DickPunkt(-12,0)
\DickPunkt(-10,0)
\DickPunkt(-8,0)
\DickPunkt(-6,0)
\DickPunkt(-4,0)
\DickPunkt(-2,0)
\DickPunkt(2,0)
\DickPunkt(4,0)
\DickPunkt(6,0)
\DickPunkt(8,0)
\DickPunkt(10,0)
\Label\u{A_0}(-2,0)
\Label\u{A_1}(-4,0)
\Label\u{A_2}(-6,0)
\Label\u{A_3}(-8,0)
\Label\u{A_4}(-10,0)
\Label\u{E_0}(2,0)
\Label\u{E_1}(4,0)
\Label\u{E_2}(6,0)
\Label\u{E_3}(8,0)
\Label\u{E_4}(10,0)
$$
\centerline{\eightpoint Figure \FD}
\endinsert

Here, the only freedom for the paths is in the strip between
the vertical lines $x=-2$ and $x=+2$. Again, except for some top paths
that have two up-steps followed by two down-steps in that strip,
what happens in the substrips to the left and to the right of
the $y$-axis is independent. We then map the configurations
in these substrips to pillars, as is illustrated in Figure~\FE\ 
for the example in Figure~\FD. The weight contribution of the
uniquely determined parts to the left of $x=-2$ and to the right
of $x=+2$ is exactly the right-hand side in (\AF) and thus
equal to the denominator in (\AI).

\midinsert
$$
\Gitter(1,8)(-2,0)
\Koordinatenachsen(1,8)(-2,0)
\SPfad(-2,-2),2222222222\endSPfad
\Pfad(-2,4),43\endPfad
\Pfad(-2,2),43\endPfad
\PfadDicke{3pt}
\Pfad(-2,0),11\endPfad
\Pfad(-2,6),11\endPfad
\DickPunkt(-2,0)
\Label\ro{\longleftrightarrow}(2,5)
\PfadDicke{1pt}
\Pfad(4,0),122222222\endPfad
\Pfad(4,0),222222221\endPfad
\Pfad(4,1),1\endPfad
\Pfad(4,3),1\endPfad
\Pfad(4,5),1\endPfad
\Pfad(4,6),1\endPfad
\Pfad(4,7),1\endPfad
\Label\ro{T_1}(4,1)
\Label\ro{T_3}(4,3)
\Label\ro{\al}(4,5)
\Label\ro{\al}(4,7)
\hbox{\hskip6cm}
\Gitter(1,8)(-2,0)
\Koordinatenachsen(1,8)(-2,0)
\SPfad(-2,-2),2222222222\endSPfad
\Pfad(-2,6),43\endPfad
\Pfad(-2,2),34\endPfad
\PfadDicke{3pt}
\Pfad(-2,0),11\endPfad
\Pfad(-2,4),11\endPfad
\DickPunkt(-2,0)
\Label\ro{\longleftrightarrow}(2,5)
\PfadDicke{1pt}
\Pfad(4,0),122222222\endPfad
\Pfad(4,0),222222221\endPfad
\Pfad(4,1),1\endPfad
\Pfad(4,2),1\endPfad
\Pfad(4,4),1\endPfad
\Pfad(4,5),1\endPfad
\Pfad(4,7),1\endPfad
\Label\ro{\be}(4,1)
\Label\ro{T_2}(4,2)
\Label\ro{T_5}(4,5)
\Label\ro{\be}(4,7)
\hskip1cm
$$
\centerline{\eightpoint Figure \FE}
\endinsert

\midinsert
$$
\Gitter(11,12)(-10,-1)
\Koordinatenachsen(11,12)(-10,-1)
\SPfad(-2,-2),22222222222222\endSPfad
\SPfad(2,-2),22222222222222\endSPfad
\Pfad(-11,0),3333333333344444444444\endPfad
\Pfad(-9,0),333333311434444444\endPfad
\Pfad(-7,0),33333431144444\endPfad
\Pfad(-5,0),3334334444\endPfad
\Pfad(-3,0),311114\endPfad
\PfadDicke{3pt}
\Pfad(-2,1),1111\endPfad
\Pfad(-2,7),11\endPfad
\Pfad(0,5),11\endPfad
\DickPunkt(-11,0)
\DickPunkt(-9,0)
\DickPunkt(-7,0)
\DickPunkt(-5,0)
\DickPunkt(-3,0)
\DickPunkt(3,0)
\DickPunkt(5,0)
\DickPunkt(7,0)
\DickPunkt(9,0)
\DickPunkt(11,0)
\Label\u{A_0}(-3,0)
\Label\u{A_1}(-5,0)
\Label\u{A_2}(-7,0)
\Label\u{A_3}(-9,0)
\Label\u{A_4}(-11,0)
\Label\u{E_0}(3,0)
\Label\u{E_1}(5,0)
\Label\u{E_2}(7,0)
\Label\u{E_3}(9,0)
\Label\u{E_4}(11,0)
$$
\centerline{\eightpoint Figure \FF}
\endinsert

\medskip
The combinatorial model for the Hankel determinant in the numerator
of (\AJ) is completely analogous to the previous one. 
Essentially, what one has to do
is to raise all paths in the previous model 
by one unit, prepend an up-step and append
a down-step to each path. We present a corresponding illustration
in Figure~\FF\ and leave the details to the
reader.\quad \quad \qed
\enddemo

\subhead 7. Afterthoughts
\endsubhead
A reader, after having seen \cite{\MuWYAA} (presented here in
(\CA) and (\CB)) and Theorem~\UA, will ask: what about the ``next''
Hankel determinant,
$$
\det\left(\al\be\ga m_{i+j}+(\al\be+\be\ga+\ga\al)m_{i+j+1}
+(\al+\be+\ga)m_{i+j+2}+m_{i+j+3}\right)_{i,j=0}^{n-1}
\quad ?
$$
Nothing prevents us from attacking it in just the same way. 
As a matter of fact, the obvious
extension of our combinatorial model from Section~6 in terms of
non-intersecting lattice paths leads rather directly to the following
result. 

\proclaim{Theorem \TD}
Let $\al$, $\be$, and $\ga$ be indeterminates. Then,
for the setup of Section~2, 
for all positive integers $n$ we have
$$\multline 
\frac {\det\left(\al\be\ga m_{i+j}+(\al\be+\be\ga+\ga\al)m_{i+j+1}
+(\al+\be+\ga)m_{i+j+2}+m_{i+j+3}\right)_{i,j=0}^{n-1}} 
{\det\left(m_{i+j}\right)_{i,j=0}^{n-1}}\\
=\sum_{0\le k\le j\le n}
 f_{j}(\al) f_{k}(\be) f_{k}( \ga) 
\left(f_{n - j}(\be)\Big\vert_{\smallmatrix s_{l} \to s_{l + j + 1} \\
        t_{l} \to t_{l + j + 1}\endsmallmatrix}\right)
  \prod_{\ell=k}^{n-1}t_{\ell}\\
+
\sum_{0\le j< k\le n}
 f_{j}(\al) f_{j}(\be) f_{k}(\ga) 
\left(f_{n - k}(\be)\Big\vert_{\smallmatrix s_{l} \to s_{l + k + 1} \\
        t_{l} \to t_{l + k + 1}\endsmallmatrix}\right)
  \prod_{\ell=j}^{n-1}t_{\ell}.
\endmultline
\tag\ZA$$

\endproclaim

\demo{Sketch of proof}
We consider the same underlying graph as in the proof of Theorem~\UA\
in Section~6, except that we do not only introduce {\it additional\/}
horizontal steps between the vertical lines $x=-1$ and $x=0$ (weighted
by~$\al$) and between the vertical lines $x=0$ and $x=+1$ (weighted
by~$\be$) but also between the vertical lines $x=+1$ and $x=+2$,
weighted by~$\ga$. The combinatorial model for the determinant in the
numerator in the left-hand side of~(\ZA) consists of families
$\Cal P=(P_0,P_1,\dots,P_{n-1})$ of non-intersecting lattice paths, where
$P_i$ runs from $A_i=(-i-1,0)$ to $E_i=(\si(i)+2,0)$,
$i=0,1,\dots,n-1$, for some permutation~$\si$ that depends on $\Cal P$.
Figure~\FG\ shows an example for $n=9$. Again, we have indicated the
additional horizontal edges by thick segments.

The combinatorial analysis of such families of non-intersecting
lattice paths leads to the following conclusions:

\roster 
\item Being non-intersecting, the paths are uniquely determined to the
left of the vertical line $x=-1$ and to the right of the vertical line
$x=+2$. Only in the strip between these two vertical lines there is (some)
freedom. 
\item In that strip, there may be some paths on top,
$P_{n-1},\dots,P_m$, which have
there an up-step followed by a horizontal step (ordinary or
additional) followed by a down-step. (In Figure~\FG, these are the two
paths $P_8$ and~$P_7$.) 
\item Below $P_m$, exactly one of the vertices of our graph
(= lattice points with non-negative ordinate) on the vertical
line $x=0$ and exactly one of the vertices of our graph on the
vertical line $x=+1$ is not taken by any of the paths.
Let these vertices be $(0,j)$ and $(1,k)$, respectively.
(In the example of Figure~\FG, 
these are the points $(0,6)$ and $(1,7)$.)
We must have $m=\max\{j,k\}+1$.
\endroster

\midinsert
$$
\Einheit.6cm
\Gitter(11,11)(-9,-1)
\Koordinatenachsen(11,11)(-9,-1)
\SPfad(-1,-2),2222222222222\endSPfad
\SPfad(1,-2),2222222222222\endSPfad
\SPfad(2,-2),2222222222222\endSPfad
\Pfad(-9,0),3333333331444444444\endPfad
\Pfad(-8,0),33333333144444444\endPfad
\Pfad(-7,0),33333334444444\endPfad
\Pfad(-6,0),333331414444\endPfad
\Pfad(-5,0),33334444\endPfad
\Pfad(-4,0),333333444444\endPfad
\Pfad(-3,0),3341344\endPfad
\Pfad(-2,0),3331444\endPfad
\Pfad(-1,0),111\endPfad
\PfadDicke{3pt}
\Pfad(-1,0),11\endPfad
\Pfad(-1,5),1\endPfad
\Pfad(1,3),1\endPfad
\Pfad(0,8),1\endPfad
\DickPunkt(-9,0)
\DickPunkt(-8,0)
\DickPunkt(-7,0)
\DickPunkt(-6,0)
\DickPunkt(-5,0)
\DickPunkt(-4,0)
\DickPunkt(-3,0)
\DickPunkt(-2,0)
\DickPunkt(-1,0)
\DickPunkt(2,0)
\DickPunkt(3,0)
\DickPunkt(4,0)
\DickPunkt(5,0)
\DickPunkt(6,0)
\DickPunkt(7,0)
\DickPunkt(8,0)
\DickPunkt(9,0)
\DickPunkt(10,0)
\Label\u{A_0}(-1,0)
\Label\u{A_1}(-2,0)
\Label\u{A_2}(-3,0)
\Label\u{A_3}(-4,0)
\Label\u{A_4}(-5,0)
\Label\u{A_5}(-6,0)
\Label\u{A_6}(-7,0)
\Label\u{A_7}(-8,0)
\Label\u{A_8}(-9,0)
\Label\u{E_0}(2,0)
\Label\u{E_1}(3,0)
\Label\u{E_2}(4,0)
\Label\u{E_3}(5,0)
\Label\u{E_4}(6,0)
\Label\u{E_5}(7,0)
\Label\u{E_6}(8,0)
\Label\u{E_7}(9,0)
\Label\u{E_8}(10,0)
$$
\centerline{\eightpoint Figure \FG}
\endinsert

\medskip
Now let first $j<k$ (as in our example in Figure~\FG).
Then, necessarily, we must have a down-step from any 
lattice point $(0,y)$ with $j<y\le k$ to $(1,y-1)$.
In the substrip between the vertical lines $x=-1$ and $x=0$
but {\it below} $y=j$, we find a configuration that is enumerated by
$f_j(\al)$. Likewise in the substrip between $x=0$ and $x=+1$
but {\it below} $y=j$, we find a configuration that is enumerated by
$f_j(\be)$, while in the substrip between $x=+1$ and $x=+2$
but {\it below} $y=k$, we find a configuration that is enumerated by
$f_k(\ga)$. Moreover, in the substrip between $x=0$ and $x=+1$
but {\it above} $y=k$, we find a configuration that is enumerated by
$f_{n-k}(\be)$, where each $s_l$ has to be replaced by $s_{l+k+1}$
and each $t_l$ by $t_{l+k+1}$. This explains the second sum on the
right-hand side of~(\ZA).

The case where $j\ge k$ is disposed of similarly and it leads
to the first sum on the right-hand side of~(\ZA).\quad \quad \qed
\enddemo

In the case where the $s_i$'s and the $t_i$'s are constant for
$i\ge1$, the right-hand side of (\ZA) can again be evaluated
in closed form in terms of Chebyshev polynomials. In order to
have a succinct notation for the result, we introduce the
umbral notation
$$
U_\al^n\equiv
U_{n}\left(\tfrac {\al+s} {2\sqrt t}\right),\quad \text{for }n\ge0.
$$

\proclaim{Corollary \TF}
Let $\al$, $\be$, and $\ga$ be indeterminates, and in the setup of
Section~2 set
$s_i\equiv s$ and $t_i\equiv t$ for $i\ge1$. Then,
for all positive integers $n$, we have
$$\multline 
\frac {\det\left(\al\be\ga m_{i+j}+(\al\be+\be\ga+\ga\al)m_{i+j+1}
+(\al+\be+\ga)m_{i+j+2}+m_{i+j+3}\right)_{i,j=0}^{n-1}} 
{\det\left(m_{i+j}\right)_{i,j=0}^{n-1}}\\
=\frac {\text{\rm Num}(\al,\be,\ga)} {(\al-\be)(\be-\ga)(\ga-\al)},
\endmultline
\tag\ZB
$$
where
$$
\multline
\text{\rm Num}(\al,\be,\ga)=
t^{(3n+3)/2}
(U_\al-U_\be)
(U_\be-U_\ga)
(U_\ga-U_\al)\\
\times\big(1
 -t^{-1/2}(s-s_0)U_\al^{-1}+t^{-1}(t-t_0)U_\al^{-2}\big)\\
\times\big(1
 -t^{-1/2}(s-s_0)U_\be^{-1}+t^{-1}(t-t_0)U_\be^{-2}\big)\\
\times\big(1
 -t^{-1/2}(s-s_0)U_\ga^{-1}+t^{-1}(t-t_0)U_\ga^{-2}\big)
U_\al^nU_\be^nU_\ga^n.
\endmultline
$$
\endproclaim

Obviously, we could now again derive many Hankel determinant 
evaluations for the
special cases that were discussed in Section~4. We refrain from
doing so here in the interest of not going overboard at
this point, and instead content ourselves with
presenting a small list of examples:
$$\align 
\det\left(M_{i+j}+3M_{i+j+1}
+3M_{i+j+2}+M_{i+j+3}\right)_{i,j=0}^{n-1}\\
&\kern-3cm
=
\frac {(n+1)(n+2)^2(n+3)(2n+3)(2n+5)} {180},
\tag\ZC\\
\det\left(-M_{i+j}+3M_{i+j+1}
-3M_{i+j+2}+M_{i+j+3}\right)_{i,j=0}^{n-1}\\
&\kern-3cm
=
[n\text{ even}](-1)^{n/2}\left(\tfrac {n} {2}+1\right)^2,
\tag\ZD\\
\det\left(-M_{i+j}-M_{i+j+1}
+M_{i+j+2}+M_{i+j+3}\right)_{i,j=0}^{n-1}\\
&\kern-3cm
=
\cases 
(-1)^{n_0}(n_0+1)^2,&\text{if }n=2n_0,\\
(-1)^{n_0-1}\binom {2n_0+2}3,&\text{if }n=2n_0-1,
\endcases
\tag\ZE\\
\det\left(-C_{i+j}+3C_{i+j+1}
-3C_{i+j+2}+C_{i+j+3}\right)_{i,j=0}^{n-1}\\
&\kern-3cm
=
\cases 
(-1)^{n_0}n_0(2n_0+1),&\text{if }n=3n_0-1,\\
(-1)^{n_0}(2n_0+1)^2,&\text{if }n=3n_0,\\
(-1)^{n_0}(n_0+1)(2n_0+1),&\text{if }n=3n_0+1,\\
\endcases
\tag\ZF\\
\det\left(-C_{i+j+1}+3C_{i+j+2}
-3C_{i+j+3}+C_{i+j+4}\right)_{i,j=0}^{n-1}\\
&\kern-3cm
=
\cases 
(-1)^{n_0-1}n_0(2n_0+1),&\text{if }n=3n_0-1,\\
(-1)^{n_0}(n_0+1)(2n_0+1),&\text{if }n=3n_0,\\
(-1)^{n_0}(2n_0+2)^2,&\text{if }n=3n_0+1.
\endcases
\tag\ZG
\endalign$$
Interestingly, by chance (?) the evaluation in (\ZC) turns out to be
exactly the same as the evaluation of the Hankel determinant
$\det(C_{i+j+4})_{i,j=0}^{n-1}$.

\medskip
In the same way as we got Corollary~\TE\ from Corollary~\TG,
the above corollary entails that the expressions in
(\ZB) satisfy a linear recurrence of order~$8$. This
generalises \cite{\DoFSAA, Theorem~4}.

\proclaim{Corollary \TI}
With the setup in Section~2, let $s_i\equiv s$ and $t_i\equiv t$ for
$i\ge1$. Then the scaled Hankel determinants 
$$\multline
t_0^{-(n-1)}t^{-\binom {n-1}2}\\
\times
\det\left(\al\be\ga m_{i+j}+(\al\be+\be\ga+\ga\al)m_{i+j+1}
+(\al+\be+\ga)m_{i+j+2}+m_{i+j+3}\right)_{i,j=0}^{n-1}
\endmultline
\tag\ZH
$$
satisfy a linear recurrence of order~$8$ with constant coefficients.
\endproclaim

\medskip
Clearly, one could continue in this manner. However, the combinatorial
analysis would become more and more intricate. 
On the other hand, it should be
noted that Corollaries~\TG\ and~\TF\ (as well as~(\CB) with
(\EFa)) follow a pattern.
On this basis, we make the following conjecture on the evaluation
of the Hankel determinant of an arbitrary linear combination of
path generating functions in the case where the $s_i$'s and the
$t_i$'s are constant for $i\ge1$. In order to state the result,
we also employ umbral notation for the path generating functions,
namely $m^n\equiv m_n$ for all~$n$.

\proclaim{Conjecture \TH}
Let $\al_1,\al_2,\dots,\al_d$  be indeterminates, and in the setup of
Section~2 set
$s_i\equiv s$ and $t_i\equiv t$ for $i\ge1$. Then,
for all positive integers $d$ and $n$, we have
$$
\frac {\det\left(
m^{i+j}\prod _{\ell=1} ^{d}(\al_\ell+m)\right)_{i,j=0}^{n-1}} 
{\det\left(m_{i+j}\right)_{i,j=0}^{n-1}}
=\frac {\text{\rm Num}(\al_1,\dots,\al_d)} 
{
\prod _{1\le i<j\le d} ^{}(\al_i-\al_j)},
\tag\ZI
$$
where
$$
\multline
\text{\rm Num}(\al_1,\dots,\al_d)=
t^{\frac {1} {2}\left(dn+\binom d2\right)}
\prod _{1\le i<j\le d} 
(U_{\al_i}-U_{\al_j})\\
\times
\prod _{i=1} ^{d}
\big(1
 -t^{-1/2}(s-s_0)U_{\al_i}^{-1}+t^{-1}(t-t_0)U_{\al_i}^{-2}\big)
U_{\al_i}^n.
\endmultline
$$
\endproclaim

Obviously (compare with the proof of Corollary~\TE), 
if true, then this formula implies that the scaled
determinants on the left-hand side of (\ZI) satisfy a linear
recurrence of order $2^d$ with constant coefficients.
This would be a vast generalisation of \cite{\DoFSAA, Conj.~5}.

\bigskip
\remark{Note}
After the first version of this article was written, the second author
discovered a general identity for Hankel determinants of linear combinations
of moments of orthogonal polynomials of which Conjecture~\TH\ is a
special case, and that this identity somehow exists in the folklore
of the theory of orthogonal polynomials but does not seem to be well-known. 
These findings will be
reported in the forthcoming article \cite{\KratCO}.
\endremark

\bigskip
\remark{Acknowledgement} 
This article was written while the second author was visiting the
Institut Mittag--Leffler during
the ``Algebraic and Enumerative Combinatorics" programme
(grant no. 2016-06596 of the Swedish Research Council) in Spring 2020,
organised by Sara Billey, Petter Br\"and\'en, Sylvie Corteel
and Svante Linusson.
\endremark


\Refs
\ref\no \AignAB\by Martin    Aigner \yr 2001 \paper Catalan and other
numbers: a recurrent theme\inbook Algebraic Combinatorics and Computer
Science\eds H.~Crapo, D.~Senato\publ Springer--Verlag\publaddr
Berlin\pages 347--390\endref 

\ref\no \CaYiAA\by Naiomi T. Cameron and Andrew C.M. Yip\paper Hankel
determinants of sums of consecutive  
Motzkin numbers\jour Linear Algebra Appl\.\vol 434 \yr2011\pages
712--722\endref 

\ref\no \CiglAS\by Johann    Cigler \yr 2002 \paper Some relations between
generalized Fibonacci and Catalan numbers\jour
Sitz.ber\. d\. \"Osterr\. Akad\. Wiss\.
Math\.-Na\-tur\-wiss\. Klasse\vol 211\pages 143--154\endref 

\ref\no \CiglAT\by Johann Cigler\paper A simple approach to some
Hankel determinants\paperinfo manuscript, {\tt ar$\chi$iv:0902.1650}\endref

\ref\no \CiglAU\by Johann Cigler\paper Hankel determinants of some
polynomial sequences\paperinfo manuscript, {\tt ar$\chi$iv:1211.0816}\endref

\ref\no \CvRIAA\by  Aleksandar Cvetkovi\'c, Predrag Rajkovi\'c and
Milo\v s Ivkovi\'c 
\paper Catalan numbers, and Hankel transform, and Fibonacci numbers
\jour J. Integer Seq\. \vol 5\yr 2002\pages Article~02.1.3, 8~pp\endref  

\ref\no \DoFSAA\by  Michael Dougherty, Christopher French, Benjamin
Saderholm and Wenyang Qian\paper
Hankel transforms of linear combinations of Catalan numbers
\jour J. Integer Seq\.\vol 14\yr 2011\pages Article~11.5.1, 20~pp
\endref

\ref\no\ElouAA\by Mohamed Elouafi\paper A unified approach for the Hankel
determinants of 
 classical combinatorial numbers\jour J. Math\. Anal\. Appl\.\vol 431\yr 2015
 \pages 1253--1274\endref

\ref\no \EuWYAA\by Sen-Peng Eu, Tsai-Lien Wong and Pei-Lan Yen
\paper Hankel determinants of sums of 
consecutive weighted Schr\"oder numbers\jour Linear Algebra Appl\.\vol
437\yr 2012\pages 2285--2299\endref

\ref\no \GaRaAF\by George    Gasper and Mizan Rahman \yr 2004 \book Basic
Hypergeometric Series\bookinfo second edition\publ Encyclopedia of
Math\. And Its Applications~96, Cambridge University Press\publaddr
Cambridge\endref 

\ref\no \GeViAA\by Ira Martin Gessel and Xavier Viennot \yr 1985 \paper
Binomial determinants, paths, and hook length formulae\jour Adv\. in
Math\. \vol 58\pages 300--321\endref 

\ref\no \GeViAB\by Ira Martin Gessel and Xavier Viennot \yr 1989 \paper
Determinants, paths, and plane partitions \paperinfo preprint,
1989\finalinfo available at {\tt
  http://www.cs.brandeis.edu/\~{}ira}\endref 

\ref\no \HaReAA\by Frank Harary and Ronald C. Read\paper
The enumeration of tree-like polyhexes\jour
Proc\. Edinburgh Math\. Soc\. (2)\vol 17\yr 1970\pages 1--13\endref 

\ref\no \IsStAB\by Mourad E.H. Ismail and Dennis Stanton \yr 1997 \paper
Classical orthogonal polynomials as moments \jour Can\. J. Math\.\vol
49 \pages 520--542\endref 

\ref\no \IsStAC\by Mourad E.H. Ismail and Dennis Stanton \yr 1998 \paper More
orthogonal polynomials as moments\inbook Mathematical Essays in Honor
of Gian-Carlo Rota\eds B.~E.~Sagan, R.~P.~Stanley\publ Progress in
Math., vol.~161, Birkh\"auser\publaddr Boston \pages 377--396\endref 

\ref\no \KratBN\by Christian    Krattenthaler \yr 1999 \paper Advanced
determinant calculus\jour S\'eminaire Lotharingien Combin\.\vol 42
\rm(``The Andrews Festschrift")\pages Article~B42q, 67~pp\endref 

\ref\no \KratBZ\by Christian    Krattenthaler \yr 2005 \paper Advanced
determinant calculus: a complement\jour Linear Algebra Appl\.\vol
411\pages 68--166\endref 

\ref\no \KratCI\by Christian Krattenthaler \yr 2010\paper Determinants of
(generalised) Catalan 
numbers\jour J. Statist\. Plann\. Inference \vol 140\pages 2260--2270\endref

\ref\no \KratCO\by Christian Krattenthaler \yr \paper Hankel determinants 
of linear combinations of moments of orthogonal polynomials, II
\jour preprint, {\tt ar$\chi$iv:2101.04225}\vol \pages \endref

\ref\no \LindAA\by Bernt    Lindstr\"om \yr 1973 \paper On the vector
representations of induced matroids\jour Bull\. London
Math\. Soc\.\vol 5\pages 85--90\endref 

\ref\no\MuWYAA\by Lili Mu, Yi Wang and Yeong-Nan Yeh
\paper Hankel determinants of linear
 combinations of consecutive Catalan-like numbers
\jour Discrete Math\. \vol  340
\yr 2017\pages 3097--3103\endref

\ref\no \OEIS\by N. J. A. Sloane\book
The On-Line Encyclopedia of Integer Sequences\publ
{\tt http://oeis.org/}\endref

\ref\no \StanBI\by Richard Peter Stanley \yr 1999 \book Enumerative
Combinatorics\bookinfo Vol.~2\publ Cambridge University Press\publaddr
Cambridge\endref 

\ref\no \VienAE\by Xavier    Viennot \yr 1983 \book Une th\'eorie
combinatoire des polyn\^omes orthogonaux g\'en\'eraux\publ
UQAM\publaddr Montr\'eal, Qu\'e\-bec\finalinfo
available at 
{\tt http://www.xavierviennot.org/xavier/polynomes\underbar{\ }orthogonaux.html}\endref 


\endRefs
\enddocument